\documentclass[10pt]{article}
\usepackage[T1]{fontenc}
\usepackage{lmodern}
\usepackage{enumitem,longtable}
\usepackage{amsfonts,amsmath,amssymb,amsthm,mathtools}
\usepackage{bm,mathrsfs}
\usepackage{geometry}
\usepackage{microtype}
\usepackage{hyperref}
\hypersetup{colorlinks=true,citecolor=blue,linkcolor=blue,urlcolor=blue}

\allowdisplaybreaks[4]
\numberwithin{equation}{section}
\newtheorem{theorem}{Theorem}[section]
\newtheorem{corollary}[theorem]{Corollary}
\newtheorem{lemma}[theorem]{Lemma}
\newtheorem{proposition}[theorem]{Proposition}
\theoremstyle{definition}
\newtheorem{definition}[theorem]{Definition}
\newtheorem{assumption}[theorem]{Assumption}
\newtheorem{example}[theorem]{Example}
\theoremstyle{remark}
\newtheorem{remark}[theorem]{Remark}

\newcommand{\E}{\mathbb{E}}
\newcommand{\Pp}{\mathbb{P}}
\newcommand{\R}{\mathbb{R}}
\newcommand{\B}{\mathcal{B}}
\newcommand{\Pred}{\mathcal{P}}
\newcommand{\Ltwozero}{\mathcal{L}_2^0}
\newcommand{\norm}[1]{\left\lVert #1\right\rVert}
\newcommand{\abs}[1]{\left\lvert #1\right\rvert}
\newcommand{\pair}[2]{\left\langle #1,#2\right\rangle}

\DeclareMathOperator{\Graph}{Graph}
\DeclareMathOperator{\diver}{div}

\title{A New Class of Stochastic Differential Variational--Hemivariational Inequalities in Hilbert Spaces with an Application\thanks{This work supported by  the National Natural Science Foundation of China (12171339, 12601630).}}
\author{Han-zhang Song$^1$, Yao-jia Zhang$^2$ and Nan-jing Huang$^1$\\{\small\it 1. Department of Mathematics, Sichuan University, Chengdu, Sichuan 610064, P.R. China}\\{\small\it 2. School of Science, Southwest Petroleum University, Chengdu 610500, P.R. China}}
\date{}

\begin{document}
\maketitle

\begin{abstract}
This paper introduces a new stochastic differential variational--hemivariational inequality (SDVHVI) in which an infinite-dimensional stochastic evolution equation is coupled, at almost every time and sample point, with a parameter-dependent variational--hemivariational inequality. The model simultaneously accommodates convex constraints, a state-dependent convex bifunction, a convex potential, and a locally Lipschitz nonconvex energy represented by the Clarke generalized directional derivative. Under uniform strong monotonicity, a quantitative smallness condition, and an integrable reference-point bound, we show that the static variational--hemivariational inequality has a unique solution. We establish a Lipschitz estimate for its dependence on the stochastic state, prove joint measurability of the solution map and predictability of its composition with a predictable state, and then reduce the coupled problem to a stochastic evolution equation with globally Lipschitz coefficients. A contraction argument in an exponentially weighted predictable space yields existence and uniqueness of a mild solution. As an application, we study a stochastic thermoelastic contact model with unilateral constraints, slip-dependent friction, fixed Tresca friction, and a nonmonotone normal-compliance law. Explicit trace and material constants give a mechanically interpretable smallness condition.
\end{abstract}

\noindent\textbf{Keywords:} stochastic differential variational--hemivariational inequality; random hemivariational inequality; predictable process; mild solution; stochastic thermoelastic contact problem.
\\

\noindent\textbf{Mathematics Subject Classification (2020):} 49J40, 49J53, 60H15, 47J20, 74M15.

\section{Introduction}

Many dynamical systems contain an evolving state together with a variable that
is determined at each instant by an equilibrium condition. The first variable
may describe, for example, temperature, concentration, accumulated damage, or
inventory, whereas the second may represent displacement, contact force, flow,
or an optimal response subject to constraints. Differential variational
inequalities (DVIs) provide a natural framework for such problems. They have
been used in differential complementarity systems, constrained dynamics,
dynamic games, electrical circuits, contact mechanics, and traffic networks;
see the systematic study of Pang and Stewart \cite{PangStewart2008}.
Related formulations of constrained and nonsmooth dynamics are developed in
\cite{Stewart2011}, while Brogliato and Tanwani \cite{BrogliatoTanwani2020}
survey the connections between differential variational inequalities,
complementarity systems, and dynamics with set-valued feedback. In a
Hilbert-space setting, a basic DVI seeks a state $x$ and an equilibrium
variable $u$ satisfying
\begin{equation}\label{eq:dvi}
 \left\{
 \begin{aligned}
  &\dot x(t)=f(t,x(t),u(t)),\qquad x(0)=x_0,\\
  &u(t)\in K,\\
  &\pair{G(t,x(t),u(t))}{v-u(t)}\geq0
  \qquad\text{for every }v\in K.
 \end{aligned}
 \right.
\end{equation}
The equilibrium relation belongs to the classical theory of variational
inequalities \cite{LionsStampacchia1967}. The variational inequality in
\eqref{eq:dvi} is equivalently written as
$0\in G(t,x(t),u(t))+N_K(u(t))$, where $N_K$ denotes the normal cone to
$K$. Thus $x$ evolves according to a differential equation, while $u$ is
selected from a constrained equilibrium whose data depend on the current
state. The coupling works in both directions: the value of $u$ affects the
evolution of $x$, and the new state $x$ changes the equilibrium that determines
$u$. This simple mechanism already permits changes of active constraints,
complementarity regimes, and contact states that cannot be represented by a
smooth ordinary differential equation alone.

Random perturbations are unavoidable in many of the problems modeled by
DVIs. Loads and material properties may fluctuate, observations may be noisy,
and unresolved effects may enter the state dynamics as random forcing.
Zhang et al. \cite{ZhangGouHuangLi2023} studied this situation
through stochastic differential variational inequalities (SDVIs). A typical
SDVI has the form
\begin{equation}\label{eq:intro-sdvi}
 \left\{
 \begin{aligned}
  &dx(t)=f(t,\omega,x(t),u(t))dt
  +g(t,\omega,x(t),u(t))dW(t),\\
  &x(0)=x_0,\\
  &u(t,\omega)\in K,\\
  &\pair{G(t,\omega,x(t),u(t))}{v-u(t)}\geq 0, \quad\forall v\in K.
 \end{aligned}
 \right.
\end{equation}
In this formulation, the equilibrium variable can influence both the drift
and the intensity of the noise, while the stochastic state changes the
instantaneous variational inequality. Existence, uniqueness, and continuous
dependence on parameters were established in \cite{ZhangGouHuangLi2023},
together with applications to pollution-abatement investment, batch
fermentation, and spatial price equilibrium. Penalty and Euler schemes were
developed in
\cite{ZhangChenHuangLi2023Penalty,ZhangChenHuangLi2023Euler}, and a
fractional model with L\'evy jumps was considered in
\cite{ZengZhangHuang2025}.

The variational inequality in \eqref{eq:intro-sdvi} is well suited to convex
energies and monotone constitutive laws. It is less suitable when the
underlying response is nonconvex or locally nonmonotone, as happens in
adhesion, delamination, material softening, and certain stick--slip laws.
Hemivariational inequalities, introduced by Panagiotopoulos and formulated
through Clarke's generalized differentiation
\cite{Clarke1990,Panagiotopoulos1993}, were developed to treat precisely such
phenomena. If $J$ is a locally Lipschitz, possibly nonconvex potential, then
its Clarke generalized directional derivative $J^\circ(z;h)$ replaces the
ordinary directional derivative of a smooth convex energy. A deterministic
differential hemivariational system of the type relevant here can be written
as
\begin{equation}\label{eq:intro-dhvi}
 \left\{
 \begin{aligned}
  &\dot x(t)=Ax(t)+f(t,x(t),u(t)),\qquad x(0)=x_0,\\
  &u(t)\in K,\\
  &\pair{G(t,x(t),u(t))-\ell(t)}{v-u(t)}
  +J^\circ(\gamma u(t);\gamma(v-u(t)))\\
  &\qquad{}+F(v)-F(u(t))\geq0,\qquad\forall v\in K.
 \end{aligned}
 \right.
\end{equation}
The generalized gradient $\partial J$ may be multivalued and need not be
monotone, so \eqref{eq:intro-dhvi} does not generally reduce to the
minimization of a convex functional. This additional flexibility is important
in models of unilateral contact with nonmonotone normal compliance,
nonconvex friction, adhesion, delamination, and semipermeability.
Deterministic partial differential hemivariational inequalities and
Banach-space differential hemivariational systems were studied in
\cite{LiuZengMotreanu2018,MigorskiZeng2018}; fractional differential
hemivariational inequalities and a contact application appear in
\cite{ZengLiuMigorski2018}. For deterministic differential
variational--hemivariational inequalities, Tang et al.
\cite{TangCenNguyenZeng2020} established existence and uniqueness of mild
solutions and convergence of penalty approximations.
Subsequent work has considered evolution equations coupled with constrained
VHVIs in contact mechanics \cite{MigorskiCaiDudek2023}, history-dependent
evolutionary VHVIs \cite{MigorskiBaiDudek2025}, and optimal control for
differential quasivariational--hemivariational systems
\cite{CaiMigorskiXiao2024}. Gap functions and global error bounds for
differential VHVIs were studied in \cite{CenMinNguyenYao2023}.
Fractional extensions include systems with variable constraint sets and
thermoelastic contact \cite{ChuChenHuangLi2024} and semigroup-based
differential VHVIs \cite{SousaAraujoFaizBenaissa2026}.

Existence and uniqueness results for stationary
variational--hemivariational inequalities, together with frictional contact
applications, were obtained in
\cite{HanMigorskiSofonea2014,MigorskiOchalSofonea2017}; finite element
convergence and error estimates were developed in
\cite{HanSofoneaDanan2018}. A systematic treatment of
variational--hemivariational inequalities and their applications is given in
\cite{SofoneaMigorski2018}. The applications of variational and
hemivariational methods extend beyond contact mechanics to interior and
boundary temperature control and to incompressible Stokes and Navier--Stokes
flows with nonsmooth slip conditions. These models, including frictional
contact with nonmonotone normal compliance, are developed in Han's monograph
\cite[Chapters~3, 4, and~8]{Han2024}.

For random hemivariational inequalities, Gwinner
\cite{Gwinner2025Random} regards the random outcome as a parameter in a
family of nonsmooth equilibrium problems. By recasting each pathwise
hemivariational inequality as an equilibrium problem, the Fan--KKM principle
and the extended-real-valued equilibrium theory developed in
\cite{Gwinner2022KKM} can be applied to obtain solvability and stability for
random nonmonotone laws. This approach provides an effective framework for
random static equilibria. It does not, however, include a stochastic state
whose evolution is driven by the equilibrium variable and which, in turn,
changes the equilibrium at every instant. In stochastic evolution
hemivariational inequalities, on the other hand, the Clarke term is usually
placed directly in a stochastic inclusion; fractional, neutral,
history-dependent, controllability, and optimal-control variants may be found,
for example, in \cite{VivekEtAl2025,DhanushVijayakumar2026}.
An impulsive fractional stochastic model is studied in
\cite{YanZhang2026}.

In physical applications, random forcing often acts together with several
distinct contact mechanisms. Thermoelastic contact provides a useful example:
sliding may generate heat at the interface, while the resulting temperature
variation produces thermal strain and changes the contact pressure and the
region in contact \cite{Barber1976,JohanssonKlarbring1993}. When mechanical
relaxation is faster than thermal diffusion, temperature evolves dynamically
while displacement is determined by quasistatic equilibrium. Fluctuating heat
supplies and applied loads introduce randomness into this two-way feedback. In
the reduced model considered here, the mechanical effect on the heat balance
is represented by a linearized coupling. The equilibrium must also accommodate
nonpenetration, a friction threshold that changes with the current tangential
displacement, and a normal response that may soften or exhibit adhesion
\cite{Needleman1990,Maugis1992}. The displacement-dependent friction threshold
leads naturally to a bifunction $\varphi(u,v)$, convex in the comparison
variable $v$, rather than to a fixed convex potential alone. A further convex
contribution can be represented by $F$, while a locally Lipschitz potential
$J$ describes the nonmonotone part of the contact law; see
\cite{SofoneaMigorski2018} for the underlying variational--hemivariational
framework. These considerations motivate an
extension of Gwinner's random HVI to a random variational--hemivariational
inequality coupled with a stochastic differential equation. We thus consider
the following stochastic differential variational--hemivariational system
(SDVHVI):
\begin{equation}\label{eq:intro-sdvhvi}
\left\{
\begin{aligned}
 &dx(t,\omega)=[Ax(t,\omega)+f(t,\omega,x(t,\omega),u(t,\omega))]dt
 +g(t,\omega,x(t,\omega),u(t,\omega))dW(t),\\
 &x(0,\omega)=x_0,\\
 &u(t,\omega)\in K,\\
 &\pair{\Phi(t,\omega,x(t,\omega),u(t,\omega))-\lambda(t,\omega)}{v-u(t,\omega)}
 +\varphi(t,\omega,u(t,\omega),v)-\varphi(t,\omega,u(t,\omega),u(t,\omega))\\
 &{}+J^\circ(\gamma u(t,\omega);\gamma(v-u(t,\omega)))+F(v)-F(u(t,\omega))\geq0,
 \quad\forall v\in K,\quad\text{a.e. }t\in[0,T],\quad\text{a.s. }\omega\in\Omega.
\end{aligned}
\right.
\end{equation}
Here $A$ generates the semigroup governing the free
evolution of the state, while $f$ and $g$ describe the deterministic and
stochastic interactions between $x$ and $u$. In the instantaneous equilibrium,
$\Phi$ is the principal operator, $\lambda$ is the random external load,
$\varphi$ and $F$ describe state-dependent and fixed convex effects,
respectively, and $J^\circ(\gamma u;\gamma(v-u))$ represents a nonconvex law
transmitted through the bounded linear operator $\gamma$. Thus, unlike the
SDVI in \cite{ZhangGouHuangLi2023}, the equilibrium in
\eqref{eq:intro-sdvhvi} contains both state-dependent variational effects and
a hemivariational term. Unlike Gwinner's random HVI, it is coupled to a
stochastic evolution equation through both the drift and the diffusion. When
$\varphi=J=F=0$, the model reduces to an infinite-dimensional SDVI; when the
noise is removed and the coefficients are deterministic, it becomes a
differential variational--hemivariational system. The stochastic thermoelastic
contact problem developed in Section~\ref{sec:application} gives a concrete
realization of the full model.

The study of \eqref{eq:intro-sdvhvi} requires more than the pathwise
solvability of its equilibrium relation. Since $u$ enters both the drift and
the diffusion, it must be treated as a stochastic process with the
measurability and integrability required by the evolution equation. In
particular, an equilibrium obtained separately for each $(t,\omega)$ must
give rise to a predictable, square-integrable process when the state $x$ is
predictable and square integrable. The dependence of the equilibrium on $x$
is therefore essential to the well-posedness of the coupled system.

Further difficulties arise from the infinite-dimensional setting. Unlike
the finite-dimensional models in
\cite{ZhangGouHuangLi2023,ZengZhangHuang2025}, the present system contains a
possibly unbounded operator $A$ on a Hilbert space. Under the assumptions
considered here, the state need not belong to $D(A)$, so the evolution
equation cannot in general be interpreted in the strong sense. We instead
seek a mild solution expressed through the semigroup generated by $A$,
following the semigroup framework in \cite{Pazy1983}.
This formulation brings in stochastic convolutions and requires tools from
infinite-dimensional stochastic analysis, including integration with respect
to a $Q$-Wiener process. At the same time, the equilibrium relation is a VHVI:
the bifunction $\varphi$ and the Clarke term prevent a direct use of the
projection characterization available for ordinary variational inequalities.
Consequently, the projection arguments used for finite-dimensional SDVIs
cannot be carried over directly to the present problem.

These issues guide both our assumptions and our approach. We combine
Gwinner's variational--hemivariational solvability theorem with quantitative
estimates and a measurable-graph argument to construct an equilibrium map
compatible with the mild formulation. A dominance condition on the principal
operator controls the possible loss of monotonicity caused by $\varphi$ and
$J$, yielding a unique frozen equilibrium $U(t,\omega,x)$ and a Lipschitz
bound for its dependence on $x$. An integrable reference-point condition
provides the corresponding growth estimate. The measurability assumptions
are chosen so that the graph of $U$ is measurable; together with uniqueness,
this ensures that $U(t,\omega,x(t,\omega))$ is predictable and square
integrable whenever $x$ is. The equilibrium can then be substituted into the
evolution equation. Semigroup estimates and the It\^o isometry for
$Q$-Wiener integrals allow us to control the resulting deterministic and
stochastic convolutions, and a contraction argument in an exponentially
weighted space establishes existence and uniqueness of the coupled mild
solution on the prescribed time interval.

The main contributions of the paper are as follows.
\begin{enumerate}[label=(\roman*)]
 \item We establish existence and uniqueness for the frozen random VHVI and
 derive a global Lipschitz estimate for its solution with respect to the
 stochastic state. This result accommodates simultaneously the bifunction
 $\varphi$, the convex functional $F$, and the nonconvex potential $J$.
 \item We prove a square-integrable growth estimate from an integrable
 residual at one reference point. We then establish joint measurability of the
 equilibrium solution map and predictability of its composition with an
 arbitrary predictable state process.
 \item We construct the coupled mild solution by reducing the SDVHVI to a
 stochastic evolution equation with predictable, globally Lipschitz
 coefficients. The exponentially weighted norm removes the need for a
 small-time restriction and yields uniqueness of both components.
 \item We apply the abstract results to a stochastic thermoelastic contact
 problem with displacement-dependent friction and nonmonotone normal
 compliance, and express the decisive dominance condition in terms of the
 mechanical coefficients of the model.
\end{enumerate}

The rest of the paper is organized as follows.
Section~\ref{sec:prelim} collects the notation and the necessary results from
semigroup theory, infinite-dimensional stochastic analysis, and nonsmooth
analysis. Section~\ref{sec:model} states the mild-solution concept and the
standing assumptions for the SDVHVI. Section~\ref{sec:coupled} analyzes the frozen random VHVI, establishes the measurability and integrability of its solution map, and proves existence and uniqueness of the coupled mild solution. Finally,
Section~\ref{sec:application} applies the abstract theory to a stochastic
thermoelastic contact model, and the last section gives the concluding remarks.

\section{Preliminaries}\label{sec:prelim}

Throughout, $H_1,H_2,H_3$ are separable real Hilbert spaces, $Z$ is a separable real Banach space, and $K\subset H_2$ is nonempty, closed, and convex. Inner products and duality pairings are both denoted by $\pair{\cdot}{\cdot}$ when no ambiguity arises.

\subsection{Notations}

Following the mapping-based presentation in \cite[Section 2]{ZengZhangHuang2025}, we collect both the model data in \eqref{eq:intro-sdvhvi} and the principal functional-analytic notation used below. Here $X$ and $Y$ denote normed spaces, and they are understood to be Hilbert spaces whenever Hilbert--Schmidt operators are considered.
\small
\renewcommand{\arraystretch}{1.18}
\begin{longtable}{@{}p{0.25\textwidth}p{0.68\textwidth}@{}}
\hline
Symbol & Meaning\\
\hline
\endfirsthead
\hline
Symbol & Meaning\\
\hline
\endhead
$X^*$ & the topological dual of $X$;\\
$\pair{x^*}{x}_{X^*\times X}$ & the duality pairing between $X^*$ and $X$; the subscripts are omitted when the spaces are clear;\\
$\norm{x}_X$ & the norm of $x$ in $X$;\\
$\mathcal L(X,Y)$ & the Banach space of bounded linear operators from $X$ to $Y$;\\
$\mathcal L_2(X,Y)$ & the Hilbert space of Hilbert--Schmidt operators from $X$ to $Y$;\\
$\B(X)$ & the Borel sigma-algebra of $X$;\\
$L^p(0,T;X)$ & the Bochner--Lebesgue space of $p$-integrable $X$-valued functions;\\
$C([0,T];X)$ & the space of continuous $X$-valued functions on $[0,T]$;\\
$\Pred_0$, $\Pred$ & the predictable sigma-algebra of the raw subfiltration, and its $dt\otimes\Pp$-completion;\\
$L^2_{\Pred}(0,T;X)$ & the space of predictable $X$-valued processes that are square-integrable with respect to $dt\otimes\Pp$;\\
$J^\circ(z;h)$, $\partial J(z)$ & the Clarke generalized directional derivative at $z$ in direction $h$, and the Clarke generalized gradient;\\
$K\subset H_2$ & the nonempty, closed, and convex constraint set for the equilibrium variable;\\
$A:D(A)\subset H_1\to H_1$ & the infinitesimal generator of the strongly continuous semigroup $\{S(t)\}_{t\geq0}$;\\
$x_0:\Omega\to H_1$ & the $\mathcal F_0^0$-measurable square-integrable initial state;\\
$x:[0,T]\times\Omega\to H_1$ & the stochastic state process;\\
$u:[0,T]\times\Omega\to K$ & the stochastic equilibrium process;\\
$Q\in\mathcal L(H_3)$ & the nonnegative, self-adjoint, trace-class covariance operator of $W$;\\
$W:[0,T]\times\Omega\to H_3$ & the $H_3$-valued $Q$-Wiener process;\\
$H_0=Q^{1/2}H_3$ & the reproducing space associated with the covariance operator $Q$;\\
$\Ltwozero=\mathcal L_2(H_0,H_1)$ & the space of admissible diffusion operators;\\
$f$ & $f:[0,T]\times\Omega\times H_1\times H_2\to H_1$, the drift coefficient;\\
$g$ & $g:[0,T]\times\Omega\times H_1\times H_2\to\Ltwozero$, the diffusion coefficient;\\
$\Phi$ & $\Phi:[0,T]\times\Omega\times H_1\times H_2\to H_2$, the principal equilibrium operator;\\
$\lambda$ & $\lambda:[0,T]\times\Omega\to H_2$, the random external load;\\
$\varphi$ & $\varphi:[0,T]\times\Omega\times K\times K\to\R$, the state-dependent convex bifunction;\\
$\gamma$ & $\gamma\in\mathcal L(H_2,Z)$, the bounded linear trace or observation operator;\\
$J$ & $J:Z\to\R$, the locally Lipschitz, generally nonconvex potential;\\
$F$ & $F:H_2\to\R$, the convex continuous potential;\\
$S(t)$, $M_T$ & the semigroup generated by $A$, and $M_T=\sup_{0\leq t\leq T}\norm{S(t)}$;\\
\hline
\end{longtable}
\normalsize

\subsection{Stochastic processes and stochastic integration}

Unless stated otherwise, the definitions and results concerning
$Q$-Wiener processes, infinite-dimensional stochastic integration, and
stochastic convolutions collected in this subsection are standard and follow
Da Prato and Zabczyk \cite{DaPratoZabczyk2014}; see also
Pr\'ev\^ot and R\"ockner \cite{PrevotRockner2007} for Hilbert-space
stochastic integration.

Let $(\Omega,\mathcal F,\{\mathcal F_t\}_{t\geq0},\Pp)$ be a complete filtered probability space satisfying the usual conditions. We assume that $\{\mathcal F_t\}$ is the usual augmentation of an uncompleted filtration $\{\mathcal F_t^0\}$, and denote the predictable sigma-algebra of the latter on $[0,T]\times\Omega$ by $\Pred_0$. Define $\Pred$ to be the $dt\otimes\Pp$-completion of $\Pred_0$. Pointwise measurability of the coefficients is imposed relative to $\Pred_0$. For a separable Hilbert space $H$, set $L^2_{\Pred}(0,T;H)=L^2([0,T]\times\Omega,\Pred,dt\otimes\Pp;H)$.

Two processes in this space are identified when they agree $dt\otimes\Pp$-almost everywhere. Every $\Pred$-measurable process with values in a separable Banach space admits a $\Pred_0$-measurable version.

Let $Q\in\mathcal L(H_3)$ be nonnegative, self-adjoint, and of trace
class. Since $H_3$ is separable, there are an orthonormal basis
$\{e_k\}_{k\geq1}$ of $H_3$ and numbers $q_k\geq0$ such that $ Qe_k=q_ke_k,\operatorname{Tr}Q=\sum_{k=1}^{\infty}q_k<\infty.$

\begin{definition}[$Q$-Wiener process]\label{def:Q-Wiener}
An $H_3$-valued process $W=\{W(t)\}_{t\geq0}$ is called a
$Q$-Wiener process with respect to the filtration
$\{\mathcal F_t\}_{t\geq0}$ if:
\begin{enumerate}[label=(\roman*)]
 \item $W(0)=0$ almost surely, $W$ is adapted, and its trajectories are
 continuous almost surely;
 \item for every $0\leq s<t$, the increment $W(t)-W(s)$ is independent of
 $\mathcal F_s$;
 \item $W(t)-W(s)$ is a centered Gaussian random variable in $H_3$ with
 covariance operator $(t-s)Q$, that is,
 \[
  \E\!\left[
   \pair{W(t)-W(s)}{h}_{H_3}
   \pair{W(t)-W(s)}{k}_{H_3}
  \right]
  =(t-s)\pair{Qh}{k}_{H_3}
 \]
 for all $h,k\in H_3$.
\end{enumerate}
\end{definition}

The driving process $W$ is assumed to be a $Q$-Wiener process with respect to
$\{\mathcal F_t^0\}$, and hence also with respect to its usual augmentation.

Equivalently, on a suitable stochastic basis, $W(t)=\sum_{k=1}^{\infty}\sqrt{q_k}\,\beta_k(t)e_k$, where $\{\beta_k\}_{k\geq1}$ are independent real-valued standard Brownian
motions. The series converges in
$L^2(\Omega;C([0,T];H_3))$ for every $T>0$. In particular, $\E\norm{W(t)}_{H_3}^2=t\operatorname{Tr}Q$.

Put $H_0=Q^{1/2}H_3$ and equip this space with the inner product
\[
 \pair{z_1}{z_2}_{H_0}
 =\pair{Q^{-1/2}z_1}{Q^{-1/2}z_2}_{H_3},
\]
where $Q^{-1/2}$ denotes the pseudo-inverse on
$(\ker Q)^\perp$. Then
$\{\sqrt{q_k}e_k:q_k>0\}$ is an orthonormal basis of $H_0$. We write $\Ltwozero=\mathcal L_2(H_0,H_1)$
and, for $G\in\Ltwozero$,
\begin{equation*}
 \norm{G}_{\Ltwozero}^2
 =\sum_{k:q_k>0}\norm{G(\sqrt{q_k}e_k)}_{H_1}^2
 =\norm{GQ^{1/2}}_{\mathcal L_2(H_3,H_1)}^2.
\end{equation*}
When $G$ is originally defined on $H_3$, the last expression is the
customary quantity $\operatorname{Tr}(GQG^*)$.

\begin{lemma}
\label{lem:ito-isometry}
Let $G:[0,T]\times\Omega\to\Ltwozero$ be predictable and satisfy
\[
 \E\int_0^T\norm{G(s)}_{\Ltwozero}^2\,ds<\infty.
\]
Then, for every $t\in[0,T]$, the stochastic integral
$\int_0^tG(s)\,dW(s)$ is a well-defined element of
$L^2(\Omega;H_1)$ and has the representation
\begin{equation*}
 \int_0^tG(s)\,dW(s)
 =\sum_{k:q_k>0}\int_0^t
   G(s)(\sqrt{q_k}e_k)\,d\beta_k(s),
\end{equation*}
where the series converges in $L^2(\Omega;H_1)$. Moreover,
\begin{equation*}
 \E\norm{\int_0^tG(s)\,dW(s)}_{H_1}^2
 =\E\int_0^t\norm{G(s)}_{\Ltwozero}^2\,ds.
\end{equation*}
\end{lemma}

Let $A:D(A)\subset H_1\to H_1$ generate a strongly continuous semigroup $\{S(t)\}_{t\geq0}$. If $b\in L^2_{\Pred}(0,T;H_1)$ and $\sigma\in L^2_{\Pred}(0,T;\Ltwozero)$, then the mild form of
$dX=(AX+b)dt+\sigma dW$ is
\begin{equation*}
 X(t)=S(t)X_0+\int_0^tS(t-s)b(s)\,ds
       +\int_0^tS(t-s)\sigma(s)\,dW(s).
\end{equation*}
On each finite interval, $M_T:=\sup_{0\leq t\leq T}\norm{S(t)}<\infty$.

The measurable-composition property below is used to evaluate random
coefficients and the equilibrium solution map along predictable processes.
It is a consequence of the measurability of product maps and compositions
in Kallenberg \cite[Lemmas~1.7--1.8]{Kallenberg2002}, applied on the
measurable space $([0,T]\times\Omega,\Pred)$.

\begin{lemma}\label{lem:predictable-composition}
Let $E$ and $Y$ be separable Banach spaces and let
$\mathcal A$ be either $\Pred_0$ or $\Pred$. If
$h:([0,T]\times\Omega)\times E\to Y$ is
$\mathcal A\otimes\B(E)$-to-$\B(Y)$ measurable and
$X:[0,T]\times\Omega\to E$ is $\mathcal A$-measurable, then
$(t,\omega)\longmapsto h(t,\omega,X(t,\omega))$ is
$\mathcal A$-measurable.
\end{lemma}

We shall also use the following well-known integral form of Gronwall's
inequality.

\begin{lemma}\label{lem:gronwall}
Let $h:[0,T]\to[0,\infty)$ be locally integrable. If constants $a,b\geq0$ satisfy
\[
 h(t)\leq a+b\int_0^t h(s)\,ds
 \quad\text{for almost every }t\in[0,T],
\]
then $h(t)\leq ae^{bt}$ for almost every $t$. In particular, if $a=0$, then $h=0$ almost everywhere. If the displayed inequality holds at every $t$, then so does the conclusion.
\end{lemma}

\subsection{Clarke nonsmooth analysis}

Let $J:Z\to\R$ be locally Lipschitz. The Clarke generalized directional derivative and generalized gradient are, respectively, defined by
\begin{align*}
 J^\circ(z;h)=\limsup_{\substack{y\to z\\\tau\downarrow0}}
 \frac{J(y+\tau h)-J(y)}{\tau}\; \mbox{  and  }\;
 \partial J(z)=\{\zeta\in Z^*:\ J^\circ(z;h)\geq\pair{\zeta}{h},\,\forall h\in Z\}.
\end{align*}
The properties of the Clarke derivative used below are collected in the next lemma; see \cite{Clarke1990,MigorskiOchalSofonea2013}.

\begin{lemma}\label{lem:clarke-properties}
Let $J:Z\to\R$ be locally Lipschitz. Then, for $z,h,h_1,h_2\in Z$ and $r\geq0$:
\begin{enumerate}[label=(\roman*)]
 \item $J^\circ(z;rh)=rJ^\circ(z;h)$;
 \item $J^\circ(z;h_1+h_2)\leq J^\circ(z;h_1)+J^\circ(z;h_2)$;
 \item $h\mapsto J^\circ(z;h)$ is finite, convex, and continuous;
 \item the support formula holds:
 \begin{equation}\label{eq:clarke-support}
  J^\circ(z;h)=\max_{\zeta\in\partial J(z)}\pair{\zeta}{h};
 \end{equation}
 \item $(z,h)\mapsto J^\circ(z;h)$ is upper semicontinuous on
 $Z\times Z$ and hence is $\B(Z)\otimes\B(Z)$-measurable;
 \item if $V$ is a Banach space and $\Gamma\in\mathcal L(V,Z)$, then
 $J\circ\Gamma:V\to\R$ is locally Lipschitz and the Clarke chain rule gives $ (J\circ\Gamma)^\circ(v;w)
  \leq J^\circ(\Gamma v;\Gamma w), \forall v,w\in V.$
 \end{enumerate}
If $J$ is Lipschitz with constant $L$ in a neighborhood of $z$, then
$\norm{\zeta}_{Z^*}\leq L$ for every $\zeta\in\partial J(z)$ and
$\abs{J^\circ(z;h)}\leq L\norm{h}_Z$.
\end{lemma}

The analysis uses the following one-sided relaxation of monotonicity:
\begin{equation}\label{eq:relaxed-monotonicity}
 J^\circ(z_1;z_2-z_1)+J^\circ(z_2;z_1-z_2)
 \leq c_J\norm{z_1-z_2}_Z^2.
\end{equation}
If $\partial J$ is relaxed monotone, namely
$\pair{\zeta_1-\zeta_2}{z_1-z_2}\geq-c_J\norm{z_1-z_2}_Z^2$, then \eqref{eq:relaxed-monotonicity} follows from \eqref{eq:clarke-support}.

\begin{lemma}\label{lem:relaxed-clarke}
Suppose that
\[
 \pair{\zeta_1-\zeta_2}{z_1-z_2}
 \geq-c_J\norm{z_1-z_2}_Z^2
 \quad
 (\zeta_i\in\partial J(z_i)).
\]
Then \eqref{eq:relaxed-monotonicity} holds for all $z_1,z_2\in Z$.
\end{lemma}

\begin{proof}
Fix $z_1,z_2\in Z$ and set $d=z_1-z_2$. If $d=0$, both directions in \eqref{eq:relaxed-monotonicity} vanish. Positive homogeneity gives $J^\circ(z_i;0)=0$, and the assertion is immediate.

Assume now that $d\neq0$. By the support representation \eqref{eq:clarke-support}, the Clarke directional derivatives are the support values of the generalized gradients. Since $\partial J(z_i)$ is weakly-star compact and evaluation at a fixed direction is weakly-star continuous, there exist $\zeta_1\in\partial J(z_1)$ and $\zeta_2\in\partial J(z_2)$ such that 
$$J^\circ(z_1;z_2-z_1)=\pair{\zeta_1}{-d},\quad
 J^\circ(z_2;z_1-z_2)=\pair{\zeta_2}{d}.$$
Adding the two equalities and rearranging the duality pairings yield
\begin{align*}
&J^\circ(z_1;z_2-z_1)+J^\circ(z_2;z_1-z_2)
=-\pair{\zeta_1-\zeta_2}{d}
 =-\pair{\zeta_1-\zeta_2}{z_1-z_2}.
\end{align*}
The assumed relaxed monotonicity of $\partial J$ gives $-\pair{\zeta_1-\zeta_2}{z_1-z_2}\leq c_J\norm{z_1-z_2}_Z^2$, which is exactly \eqref{eq:relaxed-monotonicity}.
\end{proof}

\subsection{A solvability theorem for variational--hemivariational inequalities}

We shall use the following result of Gwinner \cite[Theorem~6.2]{Gwinner2022KKM}.

\begin{theorem}[Gwinner's solvability theorem]\label{thm:gwinner-vhvi}
Let $V$ be a reflexive Banach space and let $C\subset V$ be nonempty, closed, and convex. Suppose that $\mathcal A:C\to V^*$ is hemicontinuous and strongly monotone with constant $m_{\mathcal A}>0$, that is,
\[
 \pair{\mathcal A u_1-\mathcal A u_2}{u_1-u_2}
 \geq m_{\mathcal A}\norm{u_1-u_2}_V^2
 \quad (u_1,u_2\in C).
\]
Let $\psi:C\times C\to\R$ be such that $\psi(u,\cdot)$ is convex and lower semicontinuous for every $u\in C$, and assume that, for some $\alpha_\psi\geq0$,
\begin{align*}
&\psi(u_1,v_2)-\psi(u_1,v_1)
 +\psi(u_2,v_1)-\psi(u_2,v_2)\leq
 \alpha_\psi\norm{u_1-u_2}_V\norm{v_1-v_2}_V
 \quad (u_i,v_i\in C).
\end{align*}
Let $\mathcal J:V\to\R$ be locally Lipschitz and suppose that
\[
 \mathcal J^\circ(u_1;u_2-u_1)
 +\mathcal J^\circ(u_2;u_1-u_2)
 \leq\alpha_{\mathcal J}\norm{u_1-u_2}_V^2
\]
for all $u_1,u_2\in C$. Finally, let $\mathcal F:V\to(-\infty,+\infty]$ be proper, convex, and lower semicontinuous. If
\begin{equation}\label{eq:gwinner-smallness}
 \alpha_\psi+\alpha_{\mathcal J}<m_{\mathcal A},
\end{equation}
then there exists a unique $\widehat u\in C$ such that
\begin{align}\label{eq:gwinner-vhvi}
&\pair{\mathcal A\widehat u}{v-\widehat u}
 +\psi(\widehat u,v)-\psi(\widehat u,\widehat u)
 +\mathcal J^\circ(\widehat u;v-\widehat u)
 +\mathcal F(v)-\mathcal F(\widehat u)\geq0
\end{align}
for every $v\in C$.
\end{theorem}

\section{The SDVHVI and standing assumptions}\label{sec:model}

We now formulate the standing assumptions for the SDVHVI system introduced in
\eqref{eq:intro-sdvhvi}, with the spaces and model data specified in
Subsection~2.1.

\begin{definition}\label{def:mild}
A pair of equivalence classes
$(x,u)\in L^2_{\Pred}(0,T;H_1)\times L^2_{\Pred}(0,T;H_2)$
is a \emph{mild solution} of \eqref{eq:intro-sdvhvi} if it admits
$\Pred_0$-measurable representatives, denoted again by $(x,u)$, such that
$u(t,\omega)\in K$, the variational--hemivariational inequality in
\eqref{eq:intro-sdvhvi} holds for every $v\in K$ for
$dt\otimes\Pp$-almost every $(t,\omega)$, and
\begin{align}
x(t)=S(t)x_0&+\int_0^tS(t-s)f(s,\omega,x(s),u(s))\,ds
+\int_0^tS(t-s)g(s,\omega,x(s),u(s))\,dW(s)\label{eq:mild}
\end{align}
for every $t\in[0,T]$, almost surely.
\end{definition}

The hypotheses are grouped by their analytic role.
Here $\B(K)=\{B\cap K:B\in\B(H_2)\}$ denotes the trace Borel
sigma-algebra on $K$. Unless explicitly stated otherwise, every product
measurability assertion below uses the norm-Borel sigma-algebras on the state
spaces.

\begin{assumption}\label{ass:stochastic}
The operator $A$ generates a strongly continuous semigroup $S$, and the
initial value has an $\mathcal F_0^0$-measurable square-integrable version. The mapping
$f$ is $\Pred_0\otimes\B(H_1)\otimes\B(H_2)\text{-to-}\B(H_1)$ measurable,
whereas $g$ is
$\Pred_0\otimes\B(H_1)\otimes\B(H_2)\text{-to-}\B(\Ltwozero)$ measurable.
For every $(t,\omega)$, both mappings are continuous in $(x,u)$. There
exist $L_f,L_g>0$ such that
\begin{align*}
 \norm{f(t,\omega,x_1,u_1)-f(t,\omega,x_2,u_2)}_{H_1}
 &\leq L_f(\norm{x_1-x_2}_{H_1}+\norm{u_1-u_2}_{H_2}),\\
 \norm{g(t,\omega,x_1,u_1)-g(t,\omega,x_2,u_2)}_{\Ltwozero}
 &\leq L_g(\norm{x_1-x_2}_{H_1}+\norm{u_1-u_2}_{H_2}).
\end{align*}
Furthermore,
\begin{equation}\label{eq:origin-integrability}
 \E\int_0^T\bigl(\norm{f(t,\omega,0,0)}_{H_1}^2+
 \norm{g(t,\omega,0,0)}_{\Ltwozero}^2\bigr)dt<\infty.
\end{equation}
\end{assumption}

\begin{assumption}\label{ass:phiop}
The mapping $\Phi:[0,T]\times\Omega\times H_1\times H_2\to H_2$ is
$\Pred_0\otimes\B(H_1)\otimes\B(H_2)$-to-$\B(H_2)$ measurable and
continuous in $(x,u)$ for fixed $(t,\omega)$. The map
$\lambda:[0,T]\times\Omega\to H_2$ is
$\Pred_0$-to-$\B(H_2)$ measurable. There are $m_\Phi>0$ and
$L_{\Phi x}\geq0$ such that
\begin{align}
 \pair{\Phi(t,\omega,x,u_1)-\Phi(t,\omega,x,u_2)}{u_1-u_2}
 &\geq m_\Phi\norm{u_1-u_2}_{H_2}^2,\label{eq:strong-mon}\\
 \norm{\Phi(t,\omega,x_1,u)-\Phi(t,\omega,x_2,u)}_{H_2}
 &\leq L_{\Phi x}\norm{x_1-x_2}_{H_1}.\label{eq:x-lip}
\end{align}
\end{assumption}

\begin{assumption}\label{ass:nonsmooth}
The function $\varphi:[0,T]\times\Omega\times K\times K\to\R$ is
$\Pred_0\otimes\B(K)\otimes\B(K)$-to-$\B(\R)$ measurable. For each
$(t,\omega,u)$, the function $v\mapsto\varphi(t,\omega,u,v)$ is convex and
norm-lower-semicontinuous on $K$. For a constant $c_\varphi\geq0$,
\begin{align}
&\varphi(t,\omega,u_1,v_2)-\varphi(t,\omega,u_1,v_1)
 +\varphi(t,\omega,u_2,v_1)-\varphi(t,\omega,u_2,v_2)\notag\\
&\hspace{35mm}\leq c_\varphi\norm{u_1-u_2}_{H_2}\norm{v_1-v_2}_{H_2}.
\label{eq:phi-cross}
\end{align}
The function $J:Z\to\R$ is locally Lipschitz. There are $d_J>0$ and
$c_J\geq0$ such that, for every $z\in Z$ and $\zeta\in\partial J(z)$,
\begin{equation}\label{eq:J-growth}
 \norm{\zeta}_{Z^*}\leq d_J(1+\norm{z}_Z),
\end{equation}
and, for every $z_1,z_2\in Z$,
\begin{equation}\label{eq:J-relaxed}
 J^\circ(z_1;z_2-z_1)+J^\circ(z_2;z_1-z_2)
 \leq c_J\norm{z_1-z_2}_Z^2.
\end{equation}
Finally, $F:H_2\to\R$ is convex and continuous.
\end{assumption}

For $p=(t,\omega,x)$ and $u,v\in K$, define
\begin{align}
\Psi_p(u,v):={}&\pair{\Phi(t,\omega,x,u)-\lambda(t,\omega)}{v-u}
 +\varphi(t,\omega,u,v)-\varphi(t,\omega,u,u)\notag\\
&+J^\circ(\gamma u;\gamma(v-u))+F(v)-F(u).
\label{eq:Psi}
\end{align}

Set
\begin{equation}\label{eq:delta}
 \delta=m_\Phi-c_\varphi-c_J\norm{\gamma}^2.
\end{equation}

\begin{assumption}\label{ass:reference-growth}
The dominance condition $\delta>0$ holds. There exist $u_0\in K$ and a
nonnegative $\Pred_0$-measurable function
$a:[0,T]\times\Omega\to\R$ such that
\begin{equation}\label{eq:a-integrable}
 \E\int_0^T a(t,\omega)^2dt<\infty,
\end{equation}
and, for $dt\otimes\Pp$-almost every $(t,\omega)$, the bound
\begin{equation}\label{eq:reference-bound}
 \Psi_{t,\omega,0}(u_0,v)
 \geq-a(t,\omega)\norm{v-u_0}_{H_2}
\end{equation}
holds simultaneously for every $v\in K$.
\end{assumption}

Based on the assumptions above, we can get the following proposition.

\begin{proposition}
\label{prop:Psi-regularity}
Under Assumptions~\ref{ass:phiop} and~\ref{ass:nonsmooth}, the mapping
$(t,\omega,x,u,v)\mapsto\Psi_{t,\omega,x}(u,v)$ is
$\Pred_0\otimes\B(H_1)\otimes\B(K)\otimes\B(K)\text{-to-}\B(\R)$ measurable.
Moreover, for every $p=(t,\omega,x)$ and $u\in K$, the
function $v\mapsto\Psi_p(u,v)$ is convex and weakly lower semicontinuous on
$K$.
\end{proposition}

\begin{proof}
The pairing in the first term of \eqref{eq:Psi} is measurable because
$\Phi$ and $\lambda$ have the measurability stated in
Assumption~\ref{ass:phiop} and the inner product of $H_2$ is continuous.
The two terms involving $\varphi$ are measurable by
Assumption~\ref{ass:nonsmooth}; for the diagonal term, use the measurable map
$(t,\omega,x,u,v)\mapsto(t,\omega,u,u)$. The map
$(u,v)\mapsto\bigl(\gamma u,\gamma(v-u)\bigr)$ is continuous from
$K\times K$ to $Z\times Z$, while
$(z,h)\mapsto J^\circ(z;h)$ is Borel measurable by
Lemma~\ref{lem:clarke-properties}. Finally, $F$ is continuous. These
observations prove the asserted product measurability of $\Psi$.

Fix $p$ and $u$. The operator term in \eqref{eq:Psi} is affine and continuous
in $v$. The function $v\mapsto\varphi(t,\omega,u,v)$ is convex and
norm-lower-semicontinuous. By Lemma~\ref{lem:clarke-properties}, the function
$v\mapsto J^\circ(\gamma u;\gamma(v-u))$ is convex and continuous, and $F$ is
convex and continuous. The remaining
terms do not depend on $v$. Hence $v\mapsto\Psi_p(u,v)$ is convex and
norm-lower-semicontinuous. Its epigraph is therefore a norm-closed convex
subset of $K\times\R$, and thus is weakly closed. This proves weak lower
semicontinuity.
\end{proof}

The following example shows that the assumptions above can hold
simultaneously.

\begin{example}[a scalar SDVHVI satisfying all assumptions]\label{ex:scalar}
Let $H_1=H_2=H_3=Z=\R$, $K=\R$, $Q=1$, and $\gamma=I$, and let
$A=-\nu I$ for some $\nu>0$, so that $S(t)=e^{-\nu t}I$. Let
$W$ be a real standard Brownian motion, take
$x_0\in L^2(\Omega,\mathcal F_0^0;\R)$, and let $\beta,\sigma,$ and $\lambda$
be $\Pred_0$-measurable real-valued processes satisfying
\begin{equation}\label{eq:scalar-data-integrability}
 \E\int_0^T\bigl(\abs{\beta(t,\omega)}^2
 +\abs{\sigma(t,\omega)}^2+\abs{\lambda(t,\omega)}^2\bigr)dt<\infty.
\end{equation}
For constants $b_x,b_u,c_x,c_u,\ell\in\R$, $r\in\R\setminus\{0\}$,
and $m,\eta,q>0$, define
\begin{equation}\label{eq:scalar-coefficients}
\left\{
\begin{aligned}
 f(t,\omega,x,u)&=\beta(t,\omega)+b_xx+b_uu,\\
 g(t,\omega,x,u)&=\sigma(t,\omega)+c_xx+c_uu,\\
 \Phi(t,\omega,x,u)&=mu+\ell x,\\
 \varphi(t,\omega,u,v)&=ruv,\\
 J(z)&=-\frac{\eta}{2}z^2,\\
 F(v)&=\frac q2v^2.
\end{aligned}
\right.
\end{equation}
Assume
\begin{equation}\label{eq:scalar-smallness}
 m>\abs r+\eta.
\end{equation}

The stochastic coefficients in \eqref{eq:scalar-coefficients} are globally Lipschitz,
and \eqref{eq:scalar-data-integrability} gives their required integrability at
the origin. The operator $\Phi$ is strongly monotone in $u$ with
$m_\Phi=m$ and Lipschitz in $x$ with $L_{\Phi x}=\abs\ell$. Since
$v\mapsto ruv$ is affine, it is convex and continuous, and
\begin{align*}
&\varphi(u_1,v_2)-\varphi(u_1,v_1)
 +\varphi(u_2,v_1)-\varphi(u_2,v_2)\\
&\qquad=r(u_1-u_2)(v_2-v_1)
 \leq\abs r\abs{u_1-u_2}\abs{v_1-v_2}.
\end{align*}
Thus $c_\varphi=\abs r$. Moreover, $J$ is smooth and nonconvex when
$\eta>0$, with $J^\circ(z;h)=-\eta zh$ and
$\partial J(z)=\{-\eta z\}$. Consequently,
$J^\circ(z_1;z_2-z_1)+J^\circ(z_2;z_1-z_2)=\eta\abs{z_1-z_2}^2$, so
\eqref{eq:J-growth}--\eqref{eq:J-relaxed} hold with, for example,
$d_J=\max\{1,\eta\}$ and $c_J=\eta$. The functional $F$ is convex and
continuous, and all measurability and semicontinuity requirements for $\Psi$
are immediate. By \eqref{eq:scalar-smallness},
$\delta=m-\abs r-\eta>0$.

Finally, choose $u_0=0$. From \eqref{eq:Psi} and
\eqref{eq:scalar-coefficients},
$\Psi_{t,\omega,0}(0,v)=-\lambda(t,\omega)v+\frac q2v^2
\geq-\abs{\lambda(t,\omega)}\abs v$.
Hence Assumption~\ref{ass:reference-growth} holds with
$a(t,\omega)=\abs{\lambda(t,\omega)}$. We have therefore verified
Assumptions~\ref{ass:stochastic}, \ref{ass:phiop}, \ref{ass:nonsmooth},
and~\ref{ass:reference-growth} simultaneously.
In fact, the frozen equilibrium is
$U(t,\omega,x)=(\lambda(t,\omega)-\ell x)/(m+r-\eta+q)$.
The denominator is positive because of \eqref{eq:scalar-smallness}, and this
formula confirms directly that a square-integrable state produces a
square-integrable equilibrium process.
\end{example}

\section{Well-posedness of the equilibrium feedback and the coupled mild solution}\label{sec:coupled}

With the SDVHVI and its assumptions now in place, we turn to the
well-posedness of the coupled system. A pathwise solution of the equilibrium
inequality is not sufficient for this purpose: its solution must depend on
$(t,\omega,x)$ in a way that is compatible with the predictability and
square-integrability required by the stochastic evolution equation.

We therefore first show that the frozen random VHVI defines a unique
solution map $U(t,\omega,x)$ that is Lipschitz in $x$, jointly measurable, and
of linear growth. It follows that a predictable state $x$ produces the
admissible equilibrium process
$u(t,\omega)=U(t,\omega,x(t,\omega))$. At the end of this section we give our main result about the existence and
uniqueness of the coupled mild solution.

\subsection{Well-posedness of the frozen random VHVI}\label{sec:static}

Fix $p=(t,\omega,x)$. The frozen problem is to find $u\in K$ such that
\begin{equation}\label{eq:frozen}
 \Psi_p(u,v)\geq 0,\quad \forall v\in K.
\end{equation}
We first show that $\Psi_p$ defined by \eqref{eq:Psi} has the following property.
\begin{lemma}\label{lem:strong}
Under Assumptions~\ref{ass:phiop} and~\ref{ass:nonsmooth}, for all
$u_1,u_2\in K$,
\begin{equation}\label{eq:Psi-strong}
 \Psi_p(u_1,u_2)+\Psi_p(u_2,u_1)
 \leq-\delta\norm{u_1-u_2}_{H_2}^2.
\end{equation}
\end{lemma}

\begin{proof}
Fix $u_1,u_2\in K$. Using the definition \eqref{eq:Psi} twice, with $(u,v)=(u_1,u_2)$ and $(u,v)=(u_2,u_1)$, and noting that the load terms $-\pair{\lambda}{u_2-u_1}-\pair{\lambda}{u_1-u_2}=0$ and the increments $F(u_2)-F(u_1)+F(u_1)-F(u_2)=0$ cancel, we obtain
\begin{align*}
\Psi_p(u_1,u_2)+\Psi_p(u_2,u_1)
={}&\pair{\Phi(p,u_1)}{u_2-u_1}+\pair{\Phi(p,u_2)}{u_1-u_2}
 +\varphi(p,u_1,u_2)-\varphi(p,u_1,u_1)\\
&+\varphi(p,u_2,u_1)-\varphi(p,u_2,u_2)
 +J^\circ(\gamma u_1;\gamma(u_2-u_1))
 +J^\circ(\gamma u_2;\gamma(u_1-u_2)).
\end{align*}
For the principal-operator terms, the strong monotonicity condition \eqref{eq:strong-mon} gives
\[
 \pair{\Phi(p,u_1)}{u_2-u_1}+\pair{\Phi(p,u_2)}{u_1-u_2}
 =-\pair{\Phi(p,u_1)-\Phi(p,u_2)}{u_1-u_2}
 \leq-m_\Phi\norm{u_1-u_2}_{H_2}^2.
\]
Next, taking $v_1=u_1$ and $v_2=u_2$ in the cross estimate \eqref{eq:phi-cross} yields
\[
 \varphi(p,u_1,u_2)-\varphi(p,u_1,u_1)
 +\varphi(p,u_2,u_1)-\varphi(p,u_2,u_2)
 \leq c_\varphi\norm{u_1-u_2}_{H_2}^2,
\]
while applying \eqref{eq:J-relaxed} with $z_i=\gamma u_i$ gives
\[
 J^\circ(\gamma u_1;\gamma(u_2-u_1))
 +J^\circ(\gamma u_2;\gamma(u_1-u_2))
 \leq c_J\norm{\gamma(u_1-u_2)}_Z^2
 \leq c_J\norm{\gamma}^2\norm{u_1-u_2}_{H_2}^2.
\]
Combining these three bounds with the decomposition above yields $$\Psi_p(u_1,u_2)+\Psi_p(u_2,u_1)\leq(-m_\Phi+c_\varphi+c_J\norm{\gamma}^2)\norm{u_1-u_2}_{H_2}^2=-\delta\norm{u_1-u_2}_{H_2}^2,$$
which proves \eqref{eq:Psi-strong}.
\end{proof}

\begin{theorem}\label{thm:frozen}
Under Assumptions~\ref{ass:phiop} and~\ref{ass:nonsmooth}, if
$\delta>0$, then problem \eqref{eq:frozen} has exactly one solution, denoted by
$U(t,\omega,x)$.
\end{theorem}

\begin{proof}
Fix $p=(t,\omega,x)$ and set $\overline\Phi_p(u)=\Phi(t,\omega,x,u)-\lambda(t,\omega)$. We first consider the auxiliary variational--hemivariational inequality of finding $u\in K$ such that, for every $v\in K$,
\[
 \pair{\overline\Phi_p(u)}{v-u}
 +\varphi(t,\omega,u,v)-\varphi(t,\omega,u,u)
 +(J\circ\gamma)^\circ(u;v-u)+F(v)-F(u)\geq0.
\]
By Lemma~\ref{lem:clarke-properties}, $J\circ\gamma$ is locally Lipschitz and, for all $z_1,z_2\in H_2$,
\begin{align*}
(J\circ\gamma)^\circ(z_1;z_2-z_1)
+(J\circ\gamma)^\circ(z_2;z_1-z_2)
&\leq J^\circ(\gamma z_1;\gamma(z_2-z_1))
+J^\circ(\gamma z_2;\gamma(z_1-z_2))\\
&\leq c_J\norm{\gamma(z_1-z_2)}_Z^2
\leq c_J\norm{\gamma}^2\norm{z_1-z_2}_{H_2}^2.
\end{align*}
Therefore, under Assumptions~\ref{ass:phiop} and~\ref{ass:nonsmooth}, the condition $\delta>0$ shows that all the hypotheses of Theorem~\ref{thm:gwinner-vhvi} are satisfied for this auxiliary problem. Hence there exists $\widehat u\in K$ such that, for every $v\in K$,
\[
 \pair{\overline\Phi_p(\widehat u)}{v-\widehat u}
 +\varphi(t,\omega,\widehat u,v)-\varphi(t,\omega,\widehat u,\widehat u)
 +(J\circ\gamma)^\circ(\widehat u;v-\widehat u)+F(v)-F(\widehat u)\geq0.
\]
Since $(J\circ\gamma)^\circ(\widehat u;v-\widehat u)\leq J^\circ(\gamma\widehat u;\gamma(v-\widehat u))$ by Lemma~\ref{lem:clarke-properties}, it follows that $\Psi_p(\widehat u,v)\geq0$ for every $v\in K$, and hence $\widehat u$ is a solution of \eqref{eq:frozen}.

For uniqueness, suppose that \eqref{eq:frozen} has two distinct solutions $\widehat u_1,\widehat u_2\in K$. Then $\Psi_p(\widehat u_1,\widehat u_2)\geq0$ and $\Psi_p(\widehat u_2,\widehat u_1)\geq0$, whereas Lemma~\ref{lem:strong} gives $\Psi_p(\widehat u_1,\widehat u_2)+\Psi_p(\widehat u_2,\widehat u_1)\leq-\delta\norm{\widehat u_1-\widehat u_2}_{H_2}^2<0$, which is a contradiction. Therefore, \eqref{eq:frozen} has a unique solution.
\end{proof}

\begin{proposition}\label{prop:lipschitz-U}
For fixed $(t,\omega)$ and $x_1,x_2\in H_1$,
\begin{equation}\label{eq:U-lipschitz}
 \norm{U(t,\omega,x_1)-U(t,\omega,x_2)}_{H_2}
 \leq\kappa_U\norm{x_1-x_2}_{H_1},\qquad
 \kappa_U=\frac{L_{\Phi x}}{\delta}.
\end{equation}
\end{proposition}

\begin{proof}
Let $p_i=(t,\omega,x_i)$ and $u_i=U(t,\omega,x_i)$ for $i=1,2$. Since $u_i$ solves the frozen problem corresponding to $p_i$, testing the inequality for $u_1$ with $v=u_2$ and the inequality for $u_2$ with $v=u_1$ gives $\Psi_{p_1}(u_1,u_2)\geq0$ and $\Psi_{p_2}(u_2,u_1)\geq0$. Adding the two inequalities and expanding $\Psi$ yield
\begin{align*}
0\leq{}&\pair{\Phi(t,\omega,x_1,u_1)-\Phi(t,\omega,x_2,u_2)}{u_2-u_1}
+\varphi(t,\omega,u_1,u_2)-\varphi(t,\omega,u_1,u_1)\\
&+\varphi(t,\omega,u_2,u_1)-\varphi(t,\omega,u_2,u_2)
+J^\circ(\gamma u_1;\gamma(u_2-u_1))
+J^\circ(\gamma u_2;\gamma(u_1-u_2)).
\end{align*}
To estimate the first term, insert and subtract $\Phi(t,\omega,x_1,u_2)$ and apply \eqref{eq:strong-mon}, \eqref{eq:x-lip}, and the Cauchy--Schwarz inequality to obtain
\begin{align*}
&\pair{\Phi(t,\omega,x_1,u_1)-\Phi(t,\omega,x_2,u_2)}{u_2-u_1}\\
&\quad=-\pair{\Phi(t,\omega,x_1,u_1)-\Phi(t,\omega,x_1,u_2)}{u_1-u_2}
+\pair{\Phi(t,\omega,x_1,u_2)-\Phi(t,\omega,x_2,u_2)}{u_2-u_1}\\
&\quad\leq-m_\Phi\norm{u_1-u_2}_{H_2}^2
+L_{\Phi x}\norm{x_1-x_2}_{H_1}\norm{u_1-u_2}_{H_2}.
\end{align*}
As in the proof of Lemma~\ref{lem:strong}, the cross condition \eqref{eq:phi-cross} with $v_1=u_1$ and $v_2=u_2$, together with \eqref{eq:J-relaxed} and the boundedness of $\gamma$, gives
\begin{align*}
\varphi(t,\omega,u_1,u_2)-\varphi(t,\omega,u_1,u_1)
+\varphi(t,\omega,u_2,u_1)-\varphi(t,\omega,u_2,u_2)
&\leq c_\varphi\norm{u_1-u_2}_{H_2}^2,\\
J^\circ(\gamma u_1;\gamma(u_2-u_1))
+J^\circ(\gamma u_2;\gamma(u_1-u_2))
&\leq c_J\norm{\gamma}^2\norm{u_1-u_2}_{H_2}^2.
\end{align*}
Substituting these three estimates into the expansion above and recalling \eqref{eq:delta}, we find 
$$\delta\norm{u_1-u_2}_{H_2}^2\leq L_{\Phi x}\norm{x_1-x_2}_{H_1}\norm{u_1-u_2}_{H_2}.$$
 If $u_1=u_2$, estimate \eqref{eq:U-lipschitz} is immediate; otherwise, division by $\norm{u_1-u_2}_{H_2}$ gives $\norm{u_1-u_2}_{H_2}\leq(L_{\Phi x}/\delta)\norm{x_1-x_2}_{H_1}$, as desired.
\end{proof}

\begin{proposition}\label{prop:growth-U}
Under Assumptions~\ref{ass:phiop}, \ref{ass:nonsmooth}, and
\ref{ass:reference-growth}, for $dt\otimes\Pp$-almost every $(t,\omega)$
and every $x\in H_1$, the frozen solution satisfies
\begin{equation}\label{eq:U-growth2}
 \norm{U(t,\omega,x)-u_0}_{H_2}
 \leq\frac{a(t,\omega)+L_{\Phi x}\norm{x}_{H_1}}{\delta},\qquad
 \norm{U(t,\omega,x)}_{H_2}
 \leq\norm{u_0}_{H_2}+\frac{a(t,\omega)+L_{\Phi x}\norm{x}_{H_1}}{\delta}.
\end{equation}
\end{proposition}

\begin{proof}
Fix $(t,\omega)$ outside the common exceptional set in Assumption~\ref{ass:reference-growth}, let $x\in H_1$, and put $u=U(t,\omega,x)$. Because the state variable enters $\Psi$ only through $\Phi$, adding and subtracting $\Phi(t,\omega,0,u_0)$ and using \eqref{eq:reference-bound}, the Cauchy--Schwarz inequality, and \eqref{eq:x-lip} imply
\begin{align*}
\Psi_{t,\omega,x}(u_0,u)
&=\Psi_{t,\omega,0}(u_0,u)
+\pair{\Phi(t,\omega,x,u_0)-\Phi(t,\omega,0,u_0)}{u-u_0}\\
&\geq-\bigl(a(t,\omega)+L_{\Phi x}\norm{x}_{H_1}\bigr)\norm{u-u_0}_{H_2}.
\end{align*}
On the other hand, testing the frozen problem for $u$ with $v=u_0$ gives $\Psi_{t,\omega,x}(u,u_0)\geq0$, which combined with Lemma~\ref{lem:strong} yields $\delta\norm{u-u_0}_{H_2}^2\leq-\Psi_{t,\omega,x}(u_0,u)$. Combining this with the lower bound above and dividing by $\norm{u-u_0}_{H_2}$ when $u\neq u_0$ proves the first inequality in \eqref{eq:U-growth2}; the case $u=u_0$ is immediate, and the second inequality follows from the triangle inequality.
\end{proof}

\subsection{Measurability and predictability}\label{sec:measurability}

To use $U$ as feedback in the stochastic equation, we must establish its
joint measurability in the predictable parameter and the state. We assume
that $([0,T]\times\Omega,\Pred_0)$ is a standard Borel space and apply a
measurable-graph argument. The conclusion then extends to the
$dt\otimes\Pp$-completion $\Pred$.

Choose any $u_*\in K$ and, for $n\in\mathbb N$, set
$K_n=K\cap\overline B_{H_2}(u_*,n)$. Since
$H_2$ is separable and reflexive, $K_n$ with its weak topology is compact and
metrizable. Its weak Borel sigma-algebra agrees with the trace of
$\B(K)$ on $K_n$.

We shall use the following singleton-section consequence of the
Lusin--Novikov uniformization theorem; see
\cite[Theorem~18.10]{Kechris1995}.

\begin{lemma}\label{lem:singleton-graph}
Let $(X,\mathcal A)$ and $(Y,\mathcal B(Y))$ be standard Borel spaces, and
let $G\in\mathcal A\otimes\mathcal B(Y)$. If every section $G_x=\{y\in Y:(x,y)\in G\}$, $x\in X$, is a singleton, then the unique mapping $s:X\to Y$ satisfying
$G=\Graph(s)$ is $\mathcal A$-to-$\mathcal B(Y)$ measurable.
\end{lemma}

\begin{theorem}\label{thm:measurable-U}
Assume that $([0,T]\times\Omega,\Pred_0)$ is a standard Borel space.
Under Assumptions~\ref{ass:phiop} and~\ref{ass:nonsmooth} and the dominance
condition $\delta>0$, the map $U:([0,T]\times\Omega)\times H_1\to H_2$ is $\Pred_0\otimes\B(H_1)$-to-$\B(H_2)$ measurable and therefore also
$\Pred\otimes\B(H_1)$-to-$\B(H_2)$ measurable after completion.
\end{theorem}

\begin{proof}
Let $q_n(t,\omega,x,u) =\sup_{v\in K_n}[-\Psi_{t,\omega,x}(u,v)]$ and 
$q=\sup_{n\geq1}q_n$. Then for fixed $(t,\omega,x,u)$, Proposition~\ref{prop:Psi-regularity} shows
that the function being maximized is weakly upper
semicontinuous on the compact metrizable set $K_n$. Thus, its supremum is
attained. For each $c\in\R$, the set $\{q_n>c\}$ is the
projection onto $([0,T]\times\Omega)\times H_1\times K$ of
\[
 \bigcup_{m\geq1}\{(t,\omega,x,u,v):v\in K_n,
       -\Psi_{t,\omega,x}(u,v)\geq c+1/m\}.
\]
Each set in this union is Borel by Proposition~\ref{prop:Psi-regularity}
and has compact sections in $v$.
The Arsenin--Kunugui measurable-projection theorem
\cite{Kechris1995} therefore shows that $q_n$, and hence $q$, is
$\Pred_0\otimes\B(H_1)\otimes\B(K)$-measurable.

The set
\[
 G=\{(t,\omega,x,u)\in([0,T]\times\Omega)\times H_1\times K:
       q(t,\omega,x,u)\leq0\}
\]
is Borel. Because $\bigcup_nK_n=K$, its section at $(t,\omega,x)$ consists
exactly of the solutions of the corresponding frozen inequality. Each
section is a singleton by Theorem~\ref{thm:frozen}. Applying
Lemma~\ref{lem:singleton-graph} with
$X=([0,T]\times\Omega)\times H_1$ and $Y=K$ shows directly that
$U$ is $\Pred_0\otimes\B(H_1)$-measurable. It is therefore also
$\Pred\otimes\B(H_1)$-measurable.
\end{proof}

\begin{proposition}\label{prop:predictable}
If $x:[0,T]\times\Omega\to H_1$ is $\Pred_0$-measurable, then
$u(t,\omega)=U(t,\omega,x(t,\omega))$ is $\Pred_0$-measurable. The same
composition is $\Pred$-measurable whenever $x$ is $\Pred$-measurable.
If $x\in L^2_{\Pred}(0,T;H_1)$, then $u\in L^2_{\Pred}(0,T;H_2)$.
\end{proposition}

\begin{proof}
Theorem~\ref{thm:measurable-U} and
Lemma~\ref{lem:predictable-composition}, applied with either $\Pred_0$ or
$\Pred$, give the stated measurability. If
$x\in L^2_{\Pred}(0,T;H_1)$, then
\eqref{eq:U-growth2} and \eqref{eq:a-integrable} give
\[
 \E\int_0^T\norm{u(t,\omega)}_{H_2}^2dt
 \leq3\E\int_0^T\left(
 \norm{u_0}_{H_2}^2+\frac{a(t,\omega)^2}{\delta^2}
 +\frac{L_{\Phi x}^2}{\delta^2}\norm{x(t,\omega)}_{H_1}^2
 \right)dt<\infty,
\]
and hence $u\in L^2_{\Pred}(0,T;H_2)$.
\end{proof}

\subsection{Reduction to the stochastic evolution equation and coupled well-posedness}\label{sec:reduction}

The preceding results show that the equilibrium relation defines a single-valued map $U(t,\omega,x)$ with precisely the stochastic regularity required in Definition~\ref{def:mild}. It remains to close the feedback loop: the state determines $u=U(t,\omega,x)$, while this value of $u$ enters the drift and diffusion that determine the state. We begin by recording the properties of the reduced coefficients.

\begin{proposition}\label{prop:reduced-coefficients}
Define the reduced coefficients
\begin{equation}\label{eq:reduced}
 \widehat f(t,\omega,x)=f(t,\omega,x,U(t,\omega,x)),\qquad
 \widehat g(t,\omega,x)=g(t,\omega,x,U(t,\omega,x)).
\end{equation}
Then they are predictable in $(t,\omega)$, Borel measurable in $x$, and satisfy
\[
 \norm{\widehat f(t,\omega,x_1)-\widehat f(t,\omega,x_2)}_{H_1}
 \leq\widehat L_f\norm{x_1-x_2}_{H_1},\qquad
 \norm{\widehat g(t,\omega,x_1)-\widehat g(t,\omega,x_2)}_{\Ltwozero}
 \leq\widehat L_g\norm{x_1-x_2}_{H_1},
\]
where $\widehat L_f=L_f(1+\kappa_U)$ and $\widehat L_g=L_g(1+\kappa_U)$. Moreover, \eqref{eq:origin-integrability}, \eqref{eq:a-integrable}, and \eqref{eq:U-growth2} imply
\[
 \E\int_0^T\bigl(\norm{\widehat f(t,\omega,0)}_{H_1}^2+
 \norm{\widehat g(t,\omega,0)}_{\Ltwozero}^2\bigr)dt<\infty.
\]
\end{proposition}

\begin{proof}
The joint measurability of $U$ from Theorem~\ref{thm:measurable-U}, Lemma~\ref{lem:predictable-composition}, and the Carath\'eodory assumptions on $f$ and $g$ imply the asserted measurability of the compositions in \eqref{eq:reduced}. For $x_1,x_2\in H_1$, Proposition~\ref{prop:lipschitz-U} gives
\[
 \norm{\widehat f(t,x_1)-\widehat f(t,x_2)}_{H_1}
 \leq L_f\bigl(\norm{x_1-x_2}_{H_1}+\norm{U(t,x_1)-U(t,x_2)}_{H_2}\bigr)
 \leq L_f(1+\kappa_U)\norm{x_1-x_2}_{H_1},
\]
and the analogous estimate for $\widehat g$, with $L_g$ in place of $L_f$, proves the Lipschitz bounds.

At the origin, Proposition~\ref{prop:growth-U} gives $\norm{U(t,\omega,0)}_{H_2}\leq\norm{u_0}_{H_2}+\delta^{-1}a(t,\omega)$. Consequently,
\begin{align*}
 \norm{\widehat f(t,\omega,0)}_{H_1}^2
 &\leq2\norm{f(t,\omega,0,0)}_{H_1}^2
 +4L_f^2\bigl(\norm{u_0}_{H_2}^2+\delta^{-2}a(t,\omega)^2\bigr),\\
 \norm{\widehat g(t,\omega,0)}_{\Ltwozero}^2
 &\leq2\norm{g(t,\omega,0,0)}_{\Ltwozero}^2
 +4L_g^2\bigl(\norm{u_0}_{H_2}^2+\delta^{-2}a(t,\omega)^2\bigr).
\end{align*}
Integration together with \eqref{eq:origin-integrability} and \eqref{eq:a-integrable} proves the square-integrability at the origin.
\end{proof}

We next establish the convolution estimates needed for the mild formulation,
including the existence of predictable versions.

\begin{lemma}\label{lem:convolutions}
Let $A$ generate a strongly continuous semigroup $S$ on $H_1$, and set
$M_T=\sup_{0\leq t\leq T}\norm{S(t)}$.
Let $b\in L^2_{\Pred}(0,T;H_1)$ and
$\sigma\in L^2_{\Pred}(0,T;\Ltwozero)$. Define
\[
 (\mathcal Db)(t)=\int_0^tS(t-s)b(s)\,ds,
 \qquad
 (\mathcal W\sigma)(t)=\int_0^tS(t-s)\sigma(s)\,dW(s).
\]
Then both convolutions admit predictable versions and, for $t\in[0,T]$,
\[
 \E\norm{(\mathcal Db)(t)}_{H_1}^2
 \leq M_T^2t\int_0^t\E\norm{b(s)}_{H_1}^2ds,\qquad
 \E\norm{(\mathcal W\sigma)(t)}_{H_1}^2
 \leq M_T^2\int_0^t\E\norm{\sigma(s)}_{\Ltwozero}^2ds.
\]
Consequently,
\[
 \E\int_0^T\norm{(\mathcal Db)(t)}_{H_1}^2dt
 \leq M_T^2T^2\E\int_0^T\norm{b(s)}_{H_1}^2ds,\qquad
 \E\int_0^T\norm{(\mathcal W\sigma)(t)}_{H_1}^2dt
 \leq M_T^2T\E\int_0^T\norm{\sigma(s)}_{\Ltwozero}^2ds.
\]
In particular, these versions belong to $L^2_{\Pred}(0,T;H_1)$.
\end{lemma}

\begin{proof}
Choose $\Pred_0$-measurable representatives of $b$ and $\sigma$.
Their integrability and the semigroup bound ensure that the Bochner and
stochastic integrals are well defined and admit $\mathcal F_t^0$-measurable
versions for each $t$, because $W$ is a $Q$-Wiener process for the raw
filtration. By the Cauchy--Schwarz inequality, Lemma~\ref{lem:ito-isometry}, and the Hilbert--Schmidt ideal inequality, we have 
\begin{align*}
 \E\norm{(\mathcal Db)(t)}_{H_1}^2
 &\leq t\,\E\int_0^t\norm{S(t-s)b(s)}_{H_1}^2\,ds
 \leq M_T^2t\int_0^t\E\norm{b(s)}_{H_1}^2\,ds,\\
 \E\norm{(\mathcal W\sigma)(t)}_{H_1}^2
 &=\E\int_0^t\norm{S(t-s)\sigma(s)}_{\Ltwozero}^2\,ds
 \leq M_T^2\int_0^t\E\norm{\sigma(s)}_{\Ltwozero}^2\,ds,
\end{align*}
which proves the pointwise bounds.

To obtain predictable versions, we first verify continuity in mean square.
For $t,s\in[0,T]$, put $K_t(s)=S(t-s)$ for $0\leq s<t$ and $K_t(s)=0$ for $t\leq s\leq T$.
If $t_n\to t$, strong continuity of $S$ implies
$K_{t_n}(s)y\to K_t(s)y$ for every $y\in H_1$ and every $s\neq t$.
Moreover, for $G\in\Ltwozero$ and an orthonormal basis $(e_j)$ of $H_0$, one has
\[
 \norm{[K_{t_n}(s)-K_t(s)]G}_{\Ltwozero}^2
 =\sum_j\norm{[K_{t_n}(s)-K_t(s)]Ge_j}_{H_1}^2
 \longrightarrow0
 \qquad(s\neq t),
\]
by dominated convergence for the series, since its summands are bounded by
$4M_T^2\norm{Ge_j}_{H_1}^2$. Thus the Cauchy--Schwarz inequality and the
It\^o isometry yield
\begin{align*}
 \E\norm{(\mathcal Db)(t_n)-(\mathcal Db)(t)}_{H_1}^2
 &\leq T\,\E\int_0^T
 \norm{[K_{t_n}(s)-K_t(s)]b(s)}_{H_1}^2\,ds\longrightarrow0,\\
 \E\norm{(\mathcal W\sigma)(t_n)-(\mathcal W\sigma)(t)}_{H_1}^2
 &=\E\int_0^T
 \norm{[K_{t_n}(s)-K_t(s)]\sigma(s)}_{\Ltwozero}^2\,ds\longrightarrow0.
\end{align*}
These limits follow from dominated convergence using the integrable
bounds $4M_T^2\norm{b(s)}_{H_1}^2$ and
$4M_T^2\norm{\sigma(s)}_{\Ltwozero}^2$.

Let $Y$ denote either convolution. Its continuity as a map
$[0,T]\to L^2(\Omega;H_1)$ is uniform. We can therefore choose deterministic
partitions $0=t_0^n<\cdots<t_{N_n}^n=T$ with mesh tending to zero such that
the adapted step processes defined by $Y_n(0)=0$ and $Y_n(t)=Y(t_k^n)$ for $t\in(t_k^n,t_{k+1}^n]$
satisfy
$\sup_{t\in[0,T]}\E\norm{Y_n(t)-Y(t)}_{H_1}^2\leq2^{-n}$.
Each $Y_n$ is $\Pred_0$-measurable. For every fixed $t$, summability of these
expectations implies $Y_n(t)\to Y(t)$ almost surely.
Define $\widetilde Y$ to be the pointwise limit of $Y_n$ wherever it exists
and zero elsewhere. Completeness and separability of $H_1$ ensure that the
convergence set and the resulting limit are $\Pred_0$-measurable. Hence
$\widetilde Y$ is a raw-predictable version of $Y$.

Finally, Tonelli's theorem applied to the nonnegative integrands gives
\begin{align*}
 \E\int_0^T\norm{(\mathcal Db)(t)}_{H_1}^2\,dt
 &\leq M_T^2\int_0^T\left(\int_s^T t\,dt\right)\E\norm{b(s)}_{H_1}^2\,ds
 \leq M_T^2T^2\,\E\int_0^T\norm{b(s)}_{H_1}^2\,ds,\\
 \E\int_0^T\norm{(\mathcal W\sigma)(t)}_{H_1}^2\,dt
 &\leq M_T^2\int_0^T(T-s)\E\norm{\sigma(s)}_{\Ltwozero}^2\,ds
 \leq M_T^2T\,\E\int_0^T\norm{\sigma(s)}_{\Ltwozero}^2\,ds,
\end{align*}
which proves the integrated bounds and the asserted membership in $L^2_{\Pred}(0,T;H_1)$.
\end{proof}

\begin{theorem}\label{thm:main}
Under Assumptions~\ref{ass:stochastic}, \ref{ass:phiop},
\ref{ass:nonsmooth}, and~\ref{ass:reference-growth}, and the standard-Borel
hypothesis of Subsection~\ref{sec:measurability}, system
\eqref{eq:intro-sdvhvi} has a unique mild solution $(x,u)\in L^2_{\Pred}(0,T;H_1)\times L^2_{\Pred}(0,T;H_2)$, whose equilibrium component is $u(t,\omega)=U(t,\omega,x(t,\omega))$.
\end{theorem}

\begin{proof}
Let $\mathcal X=L^2_{\Pred}(0,T;H_1)$ equipped with the equivalent norm $\norm{X}_\beta^2=\E\int_0^T e^{-\beta t}\norm{X(t)}_{H_1}^2dt$ for $\beta>0$. Define $\mathcal T$ by
\[
 (\mathcal TX)(t)=S(t)x_0
 +\int_0^tS(t-s)\widehat f(s,\omega,X(s))ds
 +\int_0^tS(t-s)\widehat g(s,\omega,X(s))dW(s).
\]
We first show that $\mathcal T$ maps $\mathcal X$ into itself.
Proposition~\ref{prop:reduced-coefficients} implies
\begin{align*}
\E\int_0^T\norm{\widehat f(s,\omega,X(s))}_{H_1}^2ds
&\leq2\widehat L_f^2\E\int_0^T\norm{X(s)}_{H_1}^2ds
+2\E\int_0^T\norm{\widehat f(s,\omega,0)}_{H_1}^2ds<\infty,\\
\E\int_0^T\norm{\widehat g(s,\omega,X(s))}_{\Ltwozero}^2ds
&\leq2\widehat L_g^2\E\int_0^T\norm{X(s)}_{H_1}^2ds
+2\E\int_0^T\norm{\widehat g(s,\omega,0)}_{\Ltwozero}^2ds<\infty.
\end{align*}
Therefore, by Lemma~\ref{lem:convolutions} and the bound $\E\int_0^T\norm{S(t)x_0}_{H_1}^2dt\leq M_T^2T\E\norm{x_0}_{H_1}^2<\infty$, using $(a+b+c)^2\leq3(a^2+b^2+c^2)$ yields
\begin{align*}
\E\int_0^T\norm{(\mathcal TX)(t)}_{H_1}^2dt
\leq{}&3M_T^2T\E\norm{x_0}_{H_1}^2
+3M_T^2T^2\E\int_0^T\norm{\widehat f(s,\omega,X(s))}_{H_1}^2ds\\
&+3M_T^2T\E\int_0^T\norm{\widehat g(s,\omega,X(s))}_{\Ltwozero}^2ds<\infty.
\end{align*}
Lemma~\ref{lem:convolutions} also provides raw-predictable versions of the two
convolutions, and $t\mapsto S(t)x_0$ has such a version by the
$\mathcal F_0^0$-measurability of $x_0$. Hence $\mathcal TX\in\mathcal X$.

Let $X,Y\in\mathcal X$. By the Cauchy--Schwarz inequality, the It\^o
isometry, and Proposition~\ref{prop:reduced-coefficients}, we have
\begin{align*}
\E\norm{(\mathcal TX)(t)-(\mathcal TY)(t)}_{H_1}^2
&\leq2\E\norm{\int_0^tS(t-s)
\bigl[\widehat f(s,\omega,X(s))-\widehat f(s,\omega,Y(s))\bigr]ds}_{H_1}^2\\
&\quad+2\E\norm{\int_0^tS(t-s)
\bigl[\widehat g(s,\omega,X(s))-\widehat g(s,\omega,Y(s))\bigr]dW(s)}_{H_1}^2\\
&\leq2M_T^2\bigl(T\widehat L_f^2+\widehat L_g^2\bigr)
\int_0^t\E\norm{X(s)-Y(s)}_{H_1}^2ds.
\end{align*}
Consequently, Fubini's theorem gives
\begin{align*}
\norm{\mathcal TX-\mathcal TY}_\beta^2
&\leq2M_T^2\bigl(T\widehat L_f^2+\widehat L_g^2\bigr)
\int_0^T\E\norm{X(s)-Y(s)}_{H_1}^2
\int_s^Te^{-\beta t}dt\,ds\\
&\leq\frac{2M_T^2(T\widehat L_f^2+\widehat L_g^2)}{\beta}
\int_0^Te^{-\beta s}\E\norm{X(s)-Y(s)}_{H_1}^2ds\\
&=\frac{2M_T^2(T\widehat L_f^2+\widehat L_g^2)}{\beta}
\norm{X-Y}_\beta^2.
\end{align*}
Choosing $\beta>2M_T^2(T\widehat L_f^2+\widehat L_g^2)$, we see that $\mathcal T$ is a contraction on $\mathcal X$. By the Banach
contraction theorem, there exists a unique fixed-point class $x_*\in\mathcal X$.

Choose a $\Pred_0$-measurable representative of $x_*$ and define $x$ by the right-hand
side of $\mathcal T x_*$, taking the predictable versions of the convolutions
from Lemma~\ref{lem:convolutions}. Then $x=x_*$ in $\mathcal X$ because
$\mathcal T x_*=x_*$. Hence $\widehat f(\cdot,x)=\widehat f(\cdot,x_*)$ and
$\widehat g(\cdot,x)=\widehat g(\cdot,x_*)$ $dt\otimes\Pp$-almost everywhere.
The deterministic integrals and, by the It\^o isometry, the stochastic
integrals are therefore unchanged at each deterministic time. It follows that
$x(t)=(\mathcal T x)(t)$ for every $t\in[0,T]$, almost surely.

Set $u(t,\omega)=U(t,\omega,x(t,\omega))$. By
Theorem~\ref{thm:frozen} and Proposition~\ref{prop:predictable}, the pair
$(x,u)$ belongs to $L^2_{\Pred}(0,T;H_1)\times L^2_{\Pred}(0,T;H_2)$.
Moreover, $u(t,\omega)\in K$ and
$\Psi_{t,\omega,x(t,\omega)}(u(t,\omega),v)\geq0$ for every $v\in K$
for $dt\otimes\Pp$-almost every $(t,\omega)$. The pointwise-in-time
fixed-point identity and \eqref{eq:reduced} give, for every $t\in[0,T]$
almost surely,
\begin{align*}
  x(t)={}&S(t)x_0
  +\int_0^tS(t-s)f(s,\omega,x(s),u(s))ds
  +\int_0^tS(t-s)g(s,\omega,x(s),u(s))dW(s).
\end{align*}
Thus $(x,u)$ is a mild solution of the SDVHVI. Finally, if $(\widetilde x,\widetilde u)$ is another mild solution, then the
uniqueness of the frozen problem gives $\widetilde u(t,\omega)=U(t,\omega,\widetilde x(t,\omega))$ $dt\otimes\Pp$-almost everywhere. Consequently, $\widetilde x$ is also a fixed point of $\mathcal T$ in $\mathcal X$. The
uniqueness of the fixed-point class yields $\widetilde x=x$ and
$\widetilde u=u$ in their respective predictable $L^2$ spaces. Their mild
identities also give $\widetilde x(t)=x(t)$ almost surely at every fixed $t$.
\end{proof}

\begin{remark}\label{rem:continuity}
If the stochastic convolution in \eqref{eq:mild} has a continuous modification for every relevant predictable square-integrable integrand, then the unique state has a continuous modification. Under the usual analytic-semigroup estimates, one obtains $x\in L^2(\Omega;C([0,T];H_1))$; see \cite{DaPratoZabczyk2014}. The equilibrium variable remains predictable and square-integrable, but need not be pathwise continuous unless the data are regular in time and $\omega$.
\end{remark}

\begin{corollary}\label{cor:initial}
Let $(x_i,u_i)$ correspond to initial data $x_{0i}$, $i=1,2$, with the same coefficients. Then there exists a constant $C_T>0$, depending only on the structural constants and $T$, such that
\[
 \sup_{0\leq t\leq T}\E\norm{x_1(t)-x_2(t)}_{H_1}^2
 +\E\int_0^T\norm{u_1(t)-u_2(t)}_{H_2}^2dt
 \leq C_T\E\norm{x_{01}-x_{02}}_{H_1}^2.
\]
\end{corollary}

\begin{proof}
Set $z(t)=x_1(t)-x_2(t)$, $w(t)=u_1(t)-u_2(t)$, and $z_0=x_{01}-x_{02}$. Since both equilibrium components solve the frozen problem corresponding to their respective states, Proposition~\ref{prop:lipschitz-U} gives the pointwise feedback estimate $\norm{w(t)}_{H_2}\leq\kappa_U\norm{z(t)}_{H_1}$ for $dt\otimes\Pp$-almost every $(t,\omega)$. Subtracting the two mild equations gives
\[
 z(t)=S(t)z_0
 +\int_0^tS(t-s)\Delta f(s)ds
 +\int_0^tS(t-s)\Delta g(s)dW(s),
\]
where $\Delta f(s)=f(s,x_1(s),u_1(s))-f(s,x_2(s),u_2(s))$ and $\Delta g(s)=g(s,x_1(s),u_1(s))-g(s,x_2(s),u_2(s))$. The Lipschitz assumptions and the pointwise bound on $w(s)$ imply $\norm{\Delta f(s)}_{H_1}\leq L_f(1+\kappa_U)\norm{z(s)}_{H_1}$ and $\norm{\Delta g(s)}_{\Ltwozero}\leq L_g(1+\kappa_U)\norm{z(s)}_{H_1}$. Using $\norm{a+b+c}^2\leq3(\norm a^2+\norm b^2+\norm c^2)$ in the mild difference equation, followed by the semigroup bound, Cauchy--Schwarz, and the It\^o isometry, we obtain
\begin{align*}
 \E\norm{z(t)}_{H_1}^2
 &\leq3M_T^2\E\norm{z_0}_{H_1}^2
 +3M_T^2t\int_0^t\E\norm{\Delta f(s)}_{H_1}^2ds
 +3M_T^2\int_0^t\E\norm{\Delta g(s)}_{\Ltwozero}^2ds\\
 &\leq3M_T^2\E\norm{z_0}_{H_1}^2
 +B_T\int_0^t\E\norm{z(s)}_{H_1}^2ds,
\end{align*}
where $B_T=3M_T^2(1+\kappa_U)^2(TL_f^2+L_g^2)$. Lemma~\ref{lem:gronwall} then yields $\E\norm{x_1(t)-x_2(t)}_{H_1}^2\leq3M_T^2e^{B_Tt}\E\norm{x_{01}-x_{02}}_{H_1}^2$ for every $t\in[0,T]$. Taking the supremum gives the first required bound, and integrating the pointwise bound for $w(t)$ gives
\[
 \E\int_0^T\norm{u_1(t)-u_2(t)}_{H_2}^2dt
 \leq\kappa_U^2\int_0^T\E\norm{z(t)}_{H_1}^2dt
 \leq3\kappa_U^2TM_T^2e^{B_TT}
 \E\norm{x_{01}-x_{02}}_{H_1}^2.
\]
Therefore the assertion holds with the explicit admissible constant $C_T=3M_T^2e^{B_TT}(1+\kappa_U^2T)$. In particular, the solution map
$x_0\mapsto(x,u)$ is globally Lipschitz from
$L^2(\Omega,\mathcal F_0;H_1)$ into
$L^\infty(0,T;L^2(\Omega;H_1))\times L^2_{\Pred}(0,T;H_2)$.
\end{proof}

\section{An application to the stochastic thermoelastic contact problem}\label{sec:application}

\subsection{Physical background and modeling scope}

Consider a deformable body occupying a bounded domain $D\subset\R^d$, $d\in\{2,3\}$, during the time interval $[0,T]$. The body may come into unilateral contact with a rigid foundation along a part $\Gamma_C$ of its boundary, as in classical contact mechanics \cite{Johnson1985}. The physical starting point is classical thermoelastic frictional contact. Sliding produces heat at the interface; the resulting nonuniform temperature creates thermal strain and changes the contact pressure, which in turn changes the mechanical and thermal response. Barber \cite{Barber1976} analyzed thermoelastic contact with interfacial frictional heating, while Johansson and Klarbring \cite{JohanssonKlarbring1993} combined linear thermoelasticity, frictional heat generation, heat transfer across the interface, and unilateral contact in a coupled model. More recent deterministic VHVI models couple thermal evolution and frictional contact \cite{MigorskiCaiDudek2023,ChuChenHuangLi2024}. These mechanisms occur in sliding contacts and energy-dissipation devices such as brakes and clutches.

The quasistatic approximation reflects a separation of time scales: elastic waves and local mechanical relaxation are assumed to equilibrate rapidly compared with the thermal evolution. Accordingly, displacement is determined at each instant by an equilibrium inequality, whereas temperature retains its time-dependent diffusion. This approximation is appropriate away from impact events and high-frequency vibration; genuinely dynamic contact would require an inertial equation for the displacement.

Two further physical observations motivate the nonsmooth boundary terms. First, friction is not always described adequately by a constant Coulomb coefficient. Laboratory-based rate-and-state laws show that frictional strength may evolve with slip history and may exhibit either strengthening or weakening behavior \cite{Ruina1983}. The dependence of the threshold $k$ on $\abs{u_\tau}$ below is a reduced quasistatic representation of this effect, rather than a complete rate-and-state law. Second, adhesive and cohesive interfaces commonly have a traction--separation response that rises to a peak and then decreases during softening or decohesion. Such interface behavior underlies cohesive-zone and adhesive-contact models \cite{Needleman1990,Maugis1992}. A locally Lipschitz, nonconvex potential and its Clarke generalized gradient can represent a multivalued or locally decreasing idealization of that response.

Random thermal forcing also has a physical interpretation. Fluctuating-hydrodynamic models of heat conductors distinguish internal thermal fluctuations from externally imposed random sources and study their effects on temperature correlations \cite{DiazGuileraRubi1986}. In the present Hilbert-space model, the $Q$-Wiener process is a coarse-grained representation of unresolved thermal and environmental fluctuations. The terms $Cx$ and $\sigma_0$ distinguish state-dependent and additive noise, while $D_0u$ allows the mechanical configuration to modulate the noise amplitude.

The model below should therefore be viewed as a mathematically tractable, linearized stochastic thermoelastic contact model. In particular, the term $Bu$ is a bounded linear approximation of configuration-dependent mechanical heat input near a reference operating state. The exact frictional power would involve contact traction and tangential velocity, and would lead to a nonlinear boundary heat flux outside the present globally Lipschitz framework. Likewise, $D_0u$ is a phenomenological first-order coupling, not a microscopic fluctuation--dissipation law. These simplifications retain the experimentally relevant feedback directions while allowing a rigorous SDVHVI well-posedness analysis.

The contact law combines a slip-dependent friction bound, a fixed Tresca term, and a possibly nonmonotone normal response. The last effect cannot in general be represented by a convex potential; it is described through a Clarke generalized directional derivative and leads naturally to a hemivariational inequality \cite{MigorskiOchalSofonea2013,SofoneaMigorski2018}. We use the following notations:
\begin{itemize}
\item $x(t,z)$ is the temperature variation at $z\in D$, and $u(t,z)$ is the displacement;
\item $\mathbb C$ is the elasticity tensor, $\varepsilon(u)$ is the infinitesimal strain, and $\alpha_T$ is the thermal-expansion coefficient;
\item $\lambda(t,\omega)$ denotes the random mechanical loading, while $r(t)$ and $\sigma_0(t)$ are the thermal source in the drift and the additive-noise source, respectively;
\item $K$ is the set of kinematically admissible displacements satisfying the unilateral contact restriction on $\Gamma_C$;
\item $\varphi$ models slip-dependent friction, $F$ models a fixed Tresca friction threshold, and $J$ is the locally Lipschitz potential of the nonmonotone normal-compliance law;
\item $W$ is a $Q$-Wiener process, $B$ describes a linearized mechanical contribution to the thermal drift, $C$ gives temperature-dependent noise, and $D_0$ makes the thermal noise intensity depend on the displacement.
\end{itemize}

The stochastic thermoelastic contact problem is to find an adapted temperature process $x$ and a predictable displacement process $u$ such that, for almost every $t\in[0,T]$ and almost surely,
\begin{equation}\label{eq:contact-system}
\left\{
\begin{aligned}
dx(t)={}&[\kappa\Delta x(t)-\rho x(t)+Bu(t)+r(t)]dt\\
&+[Cx(t)+D_0u(t)+\sigma_0(t)]dW(t),\\
x(0)={}&x_0,\\
u(t,\omega)\in{}&K,\\
0\leq{}&\int_D\mathbb C\varepsilon(u(t)): \varepsilon(v-u(t))\,dz
-\int_D\alpha_Tx(t)\,\diver(v-u(t))\,dz\\
&-\pair{\lambda(t,\omega)}{v-u(t)}
+\varphi(t,\omega,u(t),v)-\varphi(t,\omega,u(t),u(t))\\
&+J^\circ(\gamma_\nu u(t);\gamma_\nu(v-u(t)))
+F(v)-F(u(t)),\qquad v\in K,
\end{aligned}
\right.
\end{equation}

The first equation in \eqref{eq:contact-system} is a linearized stochastic energy balance. The term $\kappa\Delta x$ represents Fourier heat conduction, $-\rho x$ accounts for heat loss to the surroundings, and $r$ is an imposed volumetric heat source. The terms $Bu$ and $D_0u$ allow the current mechanical configuration to influence, respectively, the drift and the diffusion of the temperature. The multiplicative term $Cx$ models fluctuations whose intensity depends on the current temperature, whereas $\sigma_0$ represents an external random thermal input.

The last two lines of \eqref{eq:contact-system} express quasistatic mechanical equilibrium by the principle of virtual work. The first integral is the work of the internal elastic stress, the second is the work generated by thermal expansion, and $\pair{\lambda}{v-u}$ is the work of the applied body forces and boundary tractions. Membership in $K$ imposes nonpenetration. The increments of $\varphi$ and $F$ describe the two frictional contributions. Finally, $J^\circ$ represents the generalized virtual work of the normal contact traction and permits adhesive, softening, or other locally nonmonotone traction--separation laws.

Thus \eqref{eq:contact-system} has a closed two-way coupling. A temperature realization changes the elastic equilibrium through thermal stress; the resulting displacement changes both the deterministic and stochastic parts of the energy balance. The subsequent subsections specify the spaces, operators, contact potentials, and quantitative conditions under which this physical system fits the abstract SDVHVI framework.

\subsection{Geometry and spaces}

Let $D\subset\R^d$, $d\in\{2,3\}$, be a bounded Lipschitz domain with
$\partial D=\Gamma_D\cup\Gamma_N\cup\Gamma_C$, where the pieces are measurable, pairwise disjoint, and $\abs{\Gamma_D}>0$. The body is clamped on $\Gamma_D$, prescribed surface tractions act on $\Gamma_N$, and possible contact with the rigid foundation occurs on $\Gamma_C$. Let $\nu$ be the outward unit normal. For a boundary vector field $v$, write $v_\nu=v\cdot\nu$ and $v_\tau=v-v_\nu\nu$ for its normal and tangential components. Set
\[
 H_1=L^2(D),\qquad
 V=\{v\in H^1(D;\R^d):v=0\text{ on }\Gamma_D\}.
\]
The strain is $\varepsilon(v)=(\nabla v+\nabla v^\top)/2$. By Korn's inequality, the $H^1$ norm on $V$ is equivalent to the strain norm. Let
$\gamma_\nu:V\to L^2(\Gamma_C)$ and
$\gamma_\tau:V\to L^2(\Gamma_C;\R^d)$ be the normal and tangential trace maps, with
\begin{equation}\label{eq:trace}
 \norm{\gamma_\nu v}_{L^2(\Gamma_C)}\leq C_\nu\norm v_V,
 \qquad
 \norm{\gamma_\tau v}_{L^2(\Gamma_C)}\leq C_\tau\norm v_V.
\end{equation}
For a gap $g_c\in L^2(\Gamma_C)$, $g_c\geq0$, define
\begin{equation}\label{eq:contact-K}
 K=\{v\in V:v_\nu\leq g_c\text{ a.e. on }\Gamma_C\}.
\end{equation}
Then $K$ is nonempty, closed, and convex, and $0\in K$. We take $H_2=V$, using the Riesz isomorphism to identify $V^*$-valued mechanical operators with $V$-valued ones.
The inequality $v_\nu\leq g_c$ states that the normal displacement cannot pass through the initial gap separating the body and the foundation; this is the usual unilateral nonpenetration constraint \cite{KikuchiOden1988}. It permits separation when the inequality is strict and contact when equality holds.

\subsection{Stochastic temperature dynamics}

Let $A=\kappa\Delta$ on $L^2(D)$, with homogeneous Dirichlet data on a nontrivial thermal boundary part and homogeneous Neumann data on the remainder. It generates an analytic $C_0$-semigroup. Here $\kappa>0$ is the thermal diffusivity. Let
$B\in\mathcal L(V,L^2(D))$, $C\in\mathcal L(L^2(D),\Ltwozero)$, and $D_0\in\mathcal L(V,\Ltwozero)$. Consider
\begin{align}
 dx(t)={}&[\kappa\Delta x(t)-\rho x(t)+Bu(t)+r(t)]dt\notag\\
 &+[Cx(t)+D_0u(t)+\sigma_0(t)]dW(t),\qquad x(0)=x_0,
 \label{eq:heat}
\end{align}
where $\rho\geq0$, and $r$, $\sigma_0$ are $\Pred_0$-measurable with
\begin{equation}\label{eq:heat-data}
 \E\int_0^T\left(\norm{r(t)}_{L^2(D)}^2+
 \norm{\sigma_0(t)}_{\Ltwozero}^2\right)dt<\infty.
\end{equation}
Assumption~\ref{ass:stochastic} follows with, for example,
$L_f=\rho+\norm B$ and $L_g=\norm C+\norm{D_0}$.

Equation \eqref{eq:heat} is understood in the mild sense. Heat is smoothed by the semigroup generated by $\kappa\Delta$, the term $-\rho x$ models exchange with an ambient environment, and $r$ represents a distributed heat supply. The operator $B$ represents the first-order change of mechanical heat input with the configuration near a reference state. In the diffusion coefficient, $Cx$ gives temperature-dependent fluctuations, $D_0u$ makes the noise intensity sensitive to the mechanical configuration, and $\sigma_0$ represents an external random thermal input. Hence the mechanical variable affects both the mean thermal evolution and its uncertainty.

\subsection{Thermoelastic equilibrium and external loading}

Let the elasticity tensor $\mathbb C(z)$ be measurable, bounded, symmetric, and uniformly positive definite: for some $c_{\mathbb C}>0$,
\begin{equation}\label{eq:elasticity}
 \mathbb C(z)\xi:\xi\geq c_{\mathbb C}\abs\xi^2
\end{equation}
for every symmetric tensor $\xi$ and almost every $z\in D$. Let $\alpha_T\in L^\infty(D)$ be the thermal-expansion coefficient and define
\begin{equation}\label{eq:Phi-contact}
 \pair{\Phi(x,u)}{v}_{V^*\times V}
 =\int_D\mathbb C\varepsilon(u):\varepsilon(v)dz
 -\int_D\alpha_Tx\,\diver v\,dz.
\end{equation}
The associated linearized stress is
\[
 \sigma(x,u)=\mathbb C\varepsilon(u)-\alpha_TxI.
\]
Thus \eqref{eq:Phi-contact} is the virtual work $\int_D\sigma(x,u):\varepsilon(v)\,dz$. Its first part is the elastic restoring force. Its second part describes expansion or contraction caused by a temperature variation; this is the channel through which the stochastic state $x$ enters mechanical equilibrium.
If $\norm v_V\leq C_K\norm{\varepsilon(v)}_{L^2(D)}$, then
\begin{equation}\label{eq:contact-mon}
 m_\Phi=\frac{c_{\mathbb C}}{C_K^2},\qquad
 L_{\Phi x}=\norm{\alpha_T}_{L^\infty(D)}C_{\diver},
\end{equation}
where $\norm{\diver v}_{L^2(D)}\leq C_{\diver}\norm v_V$. Let the $\Pred_0$-measurable load $\lambda$ satisfy
\begin{equation*}
 \pair{\lambda(t,\omega)}{v}
 =\int_D b(t,\omega,z)\cdot v(z)\,dz
 +\int_{\Gamma_N}h_N(t,\omega,z)\cdot v(z)\,d\Gamma,
\end{equation*}
where $b$ is the random body-force density and $h_N$ is the prescribed random traction on $\Gamma_N$. The load is assumed to satisfy the integrability condition
\begin{equation}\label{eq:load-integrability}
 \E\int_0^T\norm{\lambda(t)}_{V^*}^2dt<\infty.
\end{equation}

The pairing $\pair{\lambda}{v-u}$ is the virtual work of the external loads. Random dependence allows, for example, uncertain body forces, fluctuating pressure, or randomly varying boundary tractions.

\subsection{Friction and nonmonotone normal response}

Let $k:[0,T]\times\Omega\times\Gamma_C\times[0,\infty)\to[0,\infty)$ be jointly measurable with respect to $\Pred_0$, the boundary sigma-algebra, and $\B([0,\infty))$. Suppose $k(t,\omega,\cdot,0)\in L^2(\Gamma_C)$ for every $(t,\omega)$ and that $k$ is uniformly Lipschitz in its last argument:
\begin{equation}\label{eq:k-lip}
 \abs{k(t,\omega,z,r_1)-k(t,\omega,z,r_2)}\leq L_k\abs{r_1-r_2}.
\end{equation}
Define the slip-dependent friction bifunction
\begin{equation}\label{eq:varphi-contact}
 \varphi(t,\omega,u,v)=\int_{\Gamma_C}
 k(t,\omega,z,\abs{u_\tau(z)})\abs{v_\tau(z)}d\Gamma.
\end{equation}
The coefficient $k$ is the local friction bound. Dependence on $\abs{u_\tau}$ permits a reduced description of slip strengthening when $k$ increases with slip and slip weakening when it decreases, consistent with the qualitative behavior represented by slip-dependent friction laws \cite{Ruina1983}. Random dependence on $(t,\omega)$ represents uncertain surface roughness, lubrication, or environmental changes at the interface.
For $u_1,u_2,v_1,v_2\in V$, the elementary inequality
$\abs{\abs a-\abs b}\leq\abs{a-b}$, H\"older's inequality, and \eqref{eq:trace} give
\begin{align*}
&\varphi(u_1,v_2)-\varphi(u_1,v_1)
 +\varphi(u_2,v_1)-\varphi(u_2,v_2)\\
=&\int_{\Gamma_C}
 [k(\abs{u_{1\tau}})-k(\abs{u_{2\tau}})]
 [\abs{v_{2\tau}}-\abs{v_{1\tau}}]d\Gamma\\
\leq& L_k\norm{(u_1-u_2)_\tau}_{L^2(\Gamma_C)}
 \norm{(v_1-v_2)_\tau}_{L^2(\Gamma_C)}\\
\leq& L_kC_\tau^2\norm{u_1-u_2}_V\norm{v_1-v_2}_V.
\end{align*}
Therefore \eqref{eq:phi-cross} holds with
\begin{equation}\label{eq:cphi-contact}
 c_\varphi=L_kC_\tau^2.
\end{equation}
For a fixed Tresca threshold $q\in L^\infty(\Gamma_C)$, $q\geq0$, put
\begin{equation}\label{eq:F-contact}
 F(v)=\int_{\Gamma_C}q(z)\abs{v_\tau(z)}d\Gamma.
\end{equation}
This functional is convex and continuous on $V$.
The two tangential terms have distinct roles. The functional $F$ supplies a fixed friction threshold $q$, whereas $\varphi$ permits the threshold to change with the current tangential displacement. Their sum can therefore represent a baseline resistance together with slip strengthening or slip weakening.

Let $j_\nu:\Gamma_C\times\R\to\R$ be measurable in $z$ and locally Lipschitz in its second variable, with $j_\nu(\cdot,0)\in L^1(\Gamma_C)$. Set
\begin{equation}\label{eq:J-contact}
 J(w)=\int_{\Gamma_C}j_\nu(z,w(z))d\Gamma,
 \qquad w\in L^2(\Gamma_C),
\end{equation}
and take $\gamma=\gamma_\nu$. Assume that, for almost every $z\in\Gamma_C$,
\begin{align}
 \abs\xi&\leq d_j(1+\abs r),&&\xi\in\partial j_\nu(z,r),\label{eq:jn-growth}\\
 (\xi_1-\xi_2)(r_1-r_2)&\geq-m_j\abs{r_1-r_2}^2,
 &&\xi_i\in\partial j_\nu(z,r_i).
 \label{eq:jn-relaxed}
\end{align}
Under the standard hypotheses for Clarke integral functionals, \eqref{eq:jn-growth}--\eqref{eq:jn-relaxed} imply \eqref{eq:J-growth} and
\begin{equation}\label{eq:J-contact-relaxed}
 J^\circ(w_1;w_2-w_1)+J^\circ(w_2;w_1-w_2)
 \leq m_j\norm{w_1-w_2}_{L^2(\Gamma_C)}^2.
\end{equation}
Condition \eqref{eq:jn-relaxed} allows the normal traction to have locally negative slope, as occurs in adhesive, softening, or microstructured interfaces modeled by cohesive traction--separation laws \cite{Needleman1990,Maugis1992}. The constant $m_j$ quantifies the largest admissible departure from monotonicity; it is therefore the destabilizing contribution of the normal contact law.
In local form, the generalized relation may be read, up to the chosen traction sign convention, as
\[
 -\sigma_\nu(t,\omega,z)\in\partial j_\nu(z,u_\nu(t,\omega,z))
 \quad\text{on }\Gamma_C.
\]
If $j_\nu$ were smooth and convex, this would be an ordinary monotone normal-compliance law. A nonconvex $j_\nu$ permits several admissible contact tractions at the same normal displacement and locally decreasing branches, which are characteristic of adhesion and softening.

\subsection{Weak contact problem and verification of the hypotheses}

The displacement problem is: for almost every $(t,\omega)$, find $u(t,\omega)\in K$ such that for every $v\in K$,
\begin{align}
&\int_D\mathbb C\varepsilon(u):\varepsilon(v-u)dz
 -\int_D\alpha_Tx\,\diver(v-u)dz-\pair{\lambda(t,\omega)}{v-u}\notag\\
&\quad+\varphi(t,\omega,u,v)-\varphi(t,\omega,u,u)
 +J^\circ(\gamma_\nu u;\gamma_\nu(v-u))+F(v)-F(u)\geq0.
\label{eq:contact-vhvi}
\end{align}
Equations \eqref{eq:heat} and \eqref{eq:contact-vhvi} give the detailed functional formulation of the system announced in \eqref{eq:contact-system}.

The identification with the abstract problem \eqref{eq:intro-sdvhvi} is
\begin{equation*}
\begin{gathered}
 H_1=L^2(D),\qquad H_2=V,\qquad Z=L^2(\Gamma_C),
 \qquad \gamma=\gamma_\nu,\\
 A=\kappa\Delta,\qquad
 f(t,x,u)=-\rho x+Bu+r(t),\qquad
 g(t,x,u)=Cx+D_0u+\sigma_0(t),
\end{gathered}
\end{equation*}
with $K$, $\Phi$, $\varphi$, $J$, and $F$ defined in
\eqref{eq:contact-K}, \eqref{eq:Phi-contact}, \eqref{eq:varphi-contact},
\eqref{eq:J-contact}, and \eqref{eq:F-contact}, respectively. In the weak inequality, the elastic term resists deformation, the thermal term and $\lambda$ drive the body, $K$ enforces nonpenetration, the two convex boundary terms resist tangential motion, and the Clarke term supplies the possibly nonmonotone normal reaction.

This formulation also clarifies the hierarchy of models. Removing $J^\circ$ suppresses the nonmonotone normal law. If, in addition, $k$ is independent of $u_\tau$, the mechanical relation becomes a classical variational inequality with a prescribed friction threshold. Removing the mechanical feedback operators $B$ and $D_0$ makes the thermal equation one-way coupled: temperature influences displacement, but displacement no longer affects temperature.

The abstract dominance condition becomes
\begin{equation}\label{eq:physical-smallness}
 \frac{c_{\mathbb C}}{C_K^2}>L_kC_\tau^2+m_jC_\nu^2.
\end{equation}
The left-hand side of \eqref{eq:physical-smallness} is the coercive bulk stiffness measured in the displacement norm. The term $L_kC_\tau^2$ measures the sensitivity of the tangential friction threshold to changes in slip, transferred from the boundary to the bulk by the tangential trace. Similarly, $m_jC_\nu^2$ measures the maximum negative slope of the normal-compliance law after transfer through the normal trace. Thus \eqref{eq:physical-smallness} requires the stabilizing elastic response to dominate both boundary effects. This balance ensures that two different displacement fields cannot satisfy the same instantaneous contact equilibrium; if the inequality is violated, the model may admit multiple quasistatic contact states.

To verify the reference bound, take $u_0=0$. The growth condition
\eqref{eq:jn-growth} and the trace inequality imply
\[
 \abs{J^\circ(0;\gamma_\nu v)}
 \leq c_0\norm v_V,\qquad
 c_0=d_j\abs{\Gamma_C}^{1/2}C_\nu.
\]
Moreover, $k,q\geq0$ gives
$\varphi(t,\omega,0,v)-\varphi(t,\omega,0,0)\geq0$ and
$F(v)-F(0)\geq0$. Since $\Phi(0,0)=0$, it follows that
\begin{align}
 \Psi_{t,\omega,0}(0,v)
 &=-\pair{\lambda(t,\omega)}{v}
 +\varphi(t,\omega,0,v)-\varphi(t,\omega,0,0)\notag\\
 &\quad+J^\circ(0;\gamma_\nu v)+F(v)-F(0)\notag\\
 &\geq-\bigl(\norm{\lambda(t,\omega)}_{V^*}+c_0\bigr)\norm v_V.
\label{eq:contact-reference-bound}
\end{align}
Thus \eqref{eq:reference-bound} holds with
\begin{equation}\label{eq:contact-reference-a}
 a(t,\omega)=\norm{\lambda(t,\omega)}_{V^*}+c_0.
\end{equation}
The load condition \eqref{eq:load-integrability} gives
\eqref{eq:a-integrable}. Under the stated Carath\'eodory hypotheses on $k$
and $j_\nu$, the contact data have the product measurability required in
Assumptions~\ref{ass:phiop} and~\ref{ass:nonsmooth};
Proposition~\ref{prop:Psi-regularity} then supplies the measurability and weak
lower semicontinuity of the resulting contact bifunction.

\begin{theorem}\label{thm:contact}
Assume $x_0\in L^2(\Omega,\mathcal F_0^0;L^2(D))$,
\eqref{eq:heat-data}, \eqref{eq:elasticity}, \eqref{eq:load-integrability},
\eqref{eq:k-lip}, \eqref{eq:jn-growth}--\eqref{eq:jn-relaxed}, the stated
measurability and integrability conditions on the contact data, and
\eqref{eq:physical-smallness}. Then the stochastic thermoelastic contact system
\eqref{eq:heat}, \eqref{eq:contact-vhvi} has a unique mild solution
\[
 (x,u)\in L^2_{\Pred}(0,T;L^2(D))\times L^2_{\Pred}(0,T;V).
\]
For almost every $(t,\omega)$, $u(t,\omega)$ is the unique displacement solving \eqref{eq:contact-vhvi}. If the heat stochastic convolution admits its standard continuous modification, then $x$ has continuous $L^2(D)$-valued trajectories.
\end{theorem}

\begin{proof}
The heat operator and \eqref{eq:heat-data} verify Assumption~\ref{ass:stochastic}. Relations \eqref{eq:contact-mon} verify Assumption~\ref{ass:phiop}. Equations \eqref{eq:cphi-contact}, \eqref{eq:F-contact}, and \eqref{eq:J-contact-relaxed} verify the structural parts of Assumption~\ref{ass:nonsmooth}. Condition \eqref{eq:physical-smallness} is exactly $\delta>0$, while \eqref{eq:contact-reference-bound}--\eqref{eq:contact-reference-a} verify the integrable reference bound. If the predictable parameter space satisfies the standard-Borel hypothesis of Subsection~\ref{sec:measurability}, Theorem~\ref{thm:main} applies.
\end{proof}

\begin{remark}
If $J=0$ and $k$ is independent of the current slip, the application reduces to an SDVI-type contact model. The Clarke term and the state-dependent bifunction provide essential additional modeling content: they encode nonmonotone normal response and displacement-dependent friction that cannot, in general, be represented by a single monotone variational inequality.
\end{remark}

\section{Conclusions}

We formulated an infinite-dimensional SDVHVI and established its well-posedness in a predictable $L^2$ setting. Four estimates drive the analysis: uniform strong monotonicity of $\Phi$; dominance
$m_\Phi>c_\varphi+c_J\norm{\gamma}^2$; the integrable reference-point estimate \eqref{eq:reference-bound} and its growth consequence \eqref{eq:U-growth2}; and the state Lipschitz bound \eqref{eq:U-lipschitz}. Together they reduce a coupled nonsmooth stochastic system to a stochastic evolution equation with predictable, globally Lipschitz coefficients.

The framework suggests several further directions. Natural extensions include history-dependent and fractional state equations, random constraint sets, jump noise, numerical splitting methods, and optimal control. Each requires preservation of the measurable feedback structure established here; pathwise solvability alone is insufficient when the equilibrium variable enters a stochastic integral.

\bibliographystyle{plain}
\bibliography{sdvhvi_references}

\end{document}